\documentclass[english,12pt]{smfart}
\usepackage{amsmath,amssymb,amscd,amsfonts,amsmath,mathrsfs,verbatim}
\usepackage[all]{xy}
\usepackage[mathcal]{euscript}
\usepackage{epsfig}
\usepackage{times}
\usepackage{hyperref}

\newcommand{\R}{{\mathbf R}}                   
\renewcommand{\H}{{\mathbf H}}                   
\newcommand{\C}{{\mathbf C}}                   
\newcommand{\Ct}{\widehat{\mathbf C}}          
\newcommand{\CP}{\mathbf{P}^1}                
\newcommand{\CPt}{\widehat{\mathbf{P}}^1}     
\newcommand{\CPtt}{\widehat{\widehat{\mathbf{P}}}^1}     
\newcommand{\E}{{\mathcal E}}                  
\newcommand{\F}{{\mathcal F}}                  
\newcommand{\Et}{\widehat{{\mathcal E}}}        
\newcommand{\Ett}{\widehat{\Et}}              
\newcommand{\Pt}{{\widehat{\mathtt P}} }                  
\renewcommand{\tt}{\widehat{\theta}}                   
\newcommand{\ttt}{\widehat{\tt}}                   
\newcommand{\ti}{\widehat{\infty}}                   
\newcommand{\rt}{\widehat{r}}                   
\newcommand{\nt}{\widehat{n}}                   
\newcommand{\pit}{\widehat{\pi}}                   
\newcommand{\Tt}{\widehat{T}}                       
\newcommand{\Dt}{\widehat{\mathtt D}}                   
\newcommand{\Zt}{\widehat{Z}}                   
\newcommand{\Wt}{\widehat{W}}                   
\newcommand{\Mt}{\widehat{M}}                   
\newcommand{\Xt}{\widehat{x}}                   
\newcommand{\Yt}{\widehat{y}}                   
\renewcommand{\P}{{\mathbf P}}                    
\renewcommand{\O}{{\mathcal O}}                 
\newcommand{\Om}{{\mathbf \Omega}}                 
\newcommand{\Mod}{{\mathcal M}}               
\newcommand{\Modt}{\widehat{\mathcal M}}      
\newcommand{\Nahm}{{\mathcal N}}              
\newcommand{\lra}{\longrightarrow}              
\renewcommand{\d}{\mbox{d}}                      
\newcommand{\ra}{\to}                   
\newcommand{\gl}{\mathfrak{g}}                
\newcommand{\add}[1]{\mathrm{Add}_{#1}}        
\newcommand{\del}[1]{\mathrm{Del}_{#1}}        
\newcommand{\para}{\mbox{par}-}               
\newcommand{\sB}{{\mathtt B}}
\newcommand{\sC}{{\mathtt C}}
\newcommand{\sD}{{\mathtt D}}
\newcommand{\sE}{{\mathtt E}}
\newcommand{\sF}{{\mathtt F}}
\newcommand{\sX}{{\mathtt X}}
\newcommand{\sY}{{\mathtt Y}}

\newtheorem{prop}{Proposition}[section]

\newtheorem{rk}[prop]{Remark}

\newtheorem{lem}[prop]{Lemma}
\newtheorem{defn}[prop]{Definition}
\newtheorem{cor}[prop]{Corollary}
\newtheorem{thm}[prop]{Theorem}

\newcommand{\extsh}[4]{{\mathcal E}xt_{#2}^{#1}(#3,#4)}
\newcommand{\torsh}[4]{{\mathcal T}or_{#1}^{#2}(#3,#4)}

\renewcommand{\AA}{{\mathbf A}}
\newcommand{\CC}{{\mathbf C}}

\newcommand{\PP}{{\mathbf{P}}}

\newcommand{\RR}{{\mathbf{R}}}

\newcommand{\eps}{\varepsilon}

\newcommand{\isom}{\cong}

\renewcommand{\L}{{L}}
\renewcommand{\S}{{\mathcal S}}

\newcommand{\oo}{\flat} 
\newcommand{\fbl}[2]{\operatorname{F-Bl}_{#1}#2}

\renewcommand{\sf}{\mathsf{F}}

\DeclareMathOperator{\spec}{Spec}

\DeclareMathOperator{\supp}{Supp}

\DeclareMathOperator{\coker}{coker}
\DeclareMathOperator{\im}{im}

\DeclareMathOperator{\Id}{Id}

\renewcommand{\div}[1]{(#1)}

\usepackage{ifthen}
\newcommand{\gr}[2]{\operatorname{gr}^{#1}_{#2}}

\newcommand{\floor}[1]{{\lfloor #1 \rfloor}}

\newcommand{\extgl}[4]{{\mathrm{Ext}}_{#2}^{#1}(#3,#4)}
\newcommand{\mf}{{\mathfrak m}}
\DeclareMathOperator{\depth}{depth}
\DeclareMathOperator{\dimh}{dh}
\DeclareMathOperator{\ev}{ev}
\newcommand{\manustatus}{1}

\newcommand{\manudraftorfinal}[2]{
\ifthenelse{\manustatus = 1}{#1}{#2}}

\newtheorem{assumption}{Assumption}
\newtheorem{assn}{Condition}

\renewcommand{\P}{{\mathcal P}}
\newcommand{\J}{{\mathcal J}}
\newcommand{\e}{\E} 
\newcommand{\sh}{{\mathcal H}}
\newcommand{\wtd}[1]{{#1}^\#}
\DeclareMathOperator{\rank}{rank}
\newcommand{\Q}{{\mathcal Q}}
\newcommand{\chip}{{\para\chi}}

\newcommand{\M}{M}
\newcommand{\blow}{\omega}

\newcommand{\n}{{n}}
\newcommand{\kk}{{\mathbf K}}

\DeclareMathOperator{\proj}{Proj}

\DeclareMathOperator{\sym}{Sym}

\newcommand{\sP}{{\mathtt P}}

\newcommand{\abs}{{int}}

\newcommand{\pmap}{{\boldsymbol \eta}}
\newcommand{\qmap}{{\boldsymbol \rho}}
\newcommand{\pmapz}{\pmap_{T^-}}
\newcommand{\pmapP}{\pmap_{T^+}}
\newcommand{\pmaptz}{\widehat{\pmap}_{\widehat T^-}}
\newcommand{\pmaptP}{\widehat{\pmap}_{\widehat T^+}}

\newcommand{\wpmapz}{\wtd\pmap_{T^-}}
\newcommand{\wpmapP}{\wtd\pmap_{T^+}}
\newcommand{\wpmaptz}{\wtd{\widehat \pmap}_{\widehat T^-}}

\newcommand{\qmapP}{\qmap_{\widehat T^-}}
\newcommand{\qmapPt}{\qmap_{T^-}}

\newcommand{\wqmapP}{\wtd\qmap_{\widehat T^-}}

\title[Algebraic Nahm Transform]{{Algebraic Nahm Transform for parabolic Higgs bundles on $\PP^1$}}
\author{K\"ur\c{s}at Aker \\ Feza G\"ursey Institute \\ Istanbul, Turkey \\ \and \\ Szil\'ard Szab\'o \\ Alfr\'ed R\'enyi Institute of Mathematics \\ Budapest, Hungary}
\date{\today}

\begin{document}

\maketitle
\tableofcontents
\begin{abstract}
In this paper, we give a completely algebraic description of 
Nahm transform for parabolic Higgs bundles on $\PP^1$.
\end{abstract}

\section*{Introduction}
Nahm transform is a non-linear analog of Fourier transform: 
Fix a  closed additive
subgroup $\Lambda\subset \R^4$ and its dual $\Lambda^*\subset (\R^4)^*$. 
Given a $\Lambda$-invariant 
Hermitian bundle with a unitary connection on $\R^4$ 
satisfying the
anti-selfduality (ASD) equations, and of finite energy 
on $\R^4/\Lambda$,
Nahm transform produces a $\Lambda^*$-invariant 
solution of ASD equations on
$(\R^4)^*$ of finite energy on $(\R^4)^*/\Lambda^*$.
\emph{Dimensional reduction} identifies the $\Lambda$-invariant solutions
of the ASD equations on $\R^4$ 
with the solutions of the reduced equations on $\R^{4}/\Lambda$. 
Nahm transform has been studied extensively for different $\Lambda\subset \R^4$
by various authors. M. Jardim's 
expository article \cite{Jarsur} is a good introduction to the topic
with a comprehensive list of references.

In this paper, we are interested in the case $\Lambda\isom \R^2$.
Identify the quotient with the complex line 
$\C$ and denote its dual by $\Ct$. Then, dimensional
reduction yields the equations
of a holomorphic Higgs bundle with a Hermitian-Einstein metric 
on the complex line $\C$. 

\cite{Sz} describes
Nahm transform for parabolic Higgs bundles with a Hermitian-Einstein metric
on $\CC\PP^1$ 
satisying  some semisimplicity and admissibility conditions. 
These parabolic Higgs bundles have at most regular singularities 
in points  at finite distance and an irregular (Poincar\'e
rank $1$) singularity at the infinity.
Then the Nahm tranform 
of a parabolic Higgs bundle $(\E,\theta)$ is defined by following the 
steps:
\begin{enumerate}
\item
Construct
an eigensheaf $M^\oo$ on an open subset $U$ of $\CP\times \CPt$,
\item Push $M^\oo$ by the projection 
$\widehat \pi:\CP\times \CPt \rightarrow \CPt$,
\item Choose the ``right'' extension $\widehat \E$ of $\widehat \pi_* (M^\oo)$ to $\CPt$.
\end{enumerate}

Our approach here is to {\em always} work with 
projective surfaces rather than open surfaces. The main question treated in
this work is 
how to define Nahm transform of stable parabolic Higgs bundles of degree $0$ 
{\em solely}  using elementary algebraic geometry. 
As our method is algebraic, we do not treat  the Hermitian 
metrics.

This method has several advantages: it is simpler, it allows one to
compute some explicit examples, and it can be carried
out under milder assumptions on the Higgs field than in the $L^{2}$ case. 
Although Nahm transform depends fundamentally on the 
admissibility condition (Condition \ref{assn:Main}), 
we are able to remove the assumption on semisimplicity 
of the Higgs field. Also, the assumption on the order of the 
poles can be removed. 

\enlargethispage{1cm}

Both authors are grateful to the Max Planck Institute of Bonn for the hospitality 
and excellent working conditions that allowed this work to be carried out.

\section{Outline of the Paper}\label{sec:outline}
First, in Subsection \ref{ssec:HiggsM} explain the notation and the 
notions we use and the conditions under which our results hold. 
Then, in Subsection \ref{ssec:results} we describe briefly the contents of the paper. 
\subsection{Notation}\label{ssec:HiggsM}
Let $X$ be a projective scheme over a field $\kk$.
Given a birational morphism $\blow: X' \lra X$, denote the total 
and proper transforms of a Cartier divisor $\sP$ by $\blow^*\sP$ and $\wtd\blow \sP$
respectively.

Given a global section $t$ of a line bundle $L$ on $X$, denote the vanishing locus
of $t$ by $\div{t}$.

Let 
$\sP$ be an effective Cartier divisor on $X$.
Denote 
\begin{itemize}
\item $\PP_X(\O\oplus\O(-\sP)):=\proj(\sym^\bullet (\O\oplus \O(-\sP))^\vee)$ by $Z^\sP$,
\item the structure morphism by $\pi_\sP : Z^\sP \lra X$,
\item the relative hyperplane bundle by $\O_{Z^\sP}(1)$,
\item the canonical section of $\O_{Z^\sP}(1)$  by $y_\sP$,
\item the canonical section of $\O_{Z^\sP}(1)\otimes \O(\sP)$ by $x_\sP$,
\item the automorphism acting on $Z^\sP$ by $(x_\sP,y_\sP) \mapsto (-x_\sP,y_\sP)$
by $(-1)_{Z^\sP}$.
\end{itemize}
We refer to the divisor $(y_\sP)$ as the infinity section and the 
divisor $(x_\sP)$ as the zero section of $Z^\sP$.

\begin{defn}
A {\em Higgs sheaf} $(\E,\theta)$ on $X$  (with polar divisor $\sP$) consists of 
a coherent sheaf $\E$ on $X$ and a homomorphism $\theta: \E\ra \E(\sP)$.
\end{defn}

A Higgs sheaf on a projective scheme $X$ determines a unique coherent sheaf
$M^\sP$ on the surface $Z^\sP$ so that $\dim M^\sP=\dim \E$, 
$\supp M^\sP\cap (y_\sP)=\emptyset$ and
$\pi_{\sP*}M^\sP=\E(\sP)$. The sheaf $M^\sP$ is called the {\em eigensheaf}
corresponding to the Higgs sheaf $(\E,\theta)$. The support of $M^\sP$ is the
{\em spectral scheme}. The sheaf $M^\sP$ fits into an exact sequence
$$\xymatrix{0  \ar[r] & \E  \ar[rr]^<<<<<<<<<{x_\sP
-y_\sP\theta} & & \E(\sP)\otimes \O_{Z^\sP}(1)  \ar[r] & M^\sP  \ar[r] & 0}.$$
Let $$\pi^H(\E,\theta):=M^\sP.$$
 
Conversely, let the Higgs bundle $\theta:\E \ra \E(\sP)$ be the push-forward of the following 
sequence by
$\pi_\sP$
$$ \Id_M\otimes x_\sP: M(-\sP) \ra M\otimes \O_{Z^\sP}(1).$$ 
Denote
$(\E,\theta)$
by $\pi_H(M,x_\sP)$, or simply by $\pi_H(M)$. It is clear that $\pi^H$ and $\pi_H$ are 
quasi-inverses. 
The correspondence extends to parabolic objects (see Definitions \ref{def:parsh} and \ref{defn:Higgssheaf}) 
in a straightforward manner. We coin this construction as {\em the standard construction} and the resulting
objects as the standard eigensheaf, the standard spectral cover etc.

Given a parabolic Higgs bundle $(\E_\bullet,\theta)$ on $\PP^1$, define 
$\F:=\ker(\E\lra \coker\theta(-\sP)_\sP)$. The results of Sections \ref{sec:qi} and \ref{sec:moduli} 
hold under the following admissibility condition for the parabolic structure: 
\begin{assn} \label{assn:Main}
For any polar point 
$p\in \sP$, the Higgs bundle $(\E,\theta)$ satisfies  one of the following
conditions: 
\begin{itemize}
\item  $\alpha_0(\E_p)>0$ and $\E_p=\F_p$, or
\item  $\alpha_0(\E_p)=0$ and $F_1\E_p=\im(\F_p \to \E_p)$.
\end{itemize}
\end{assn}

Let $M^\sP$ and $N^\sP$ be the eigensheaves corresponding to $\E$ and 
$\F$ on $Z^\sP$. Recall that $N^\sP:=\ker(M^\sP\lra M^\sP_{T^+})$, where 
$T^+=\pi^*(\sP)\cap (x_\sP)$.
Then Condition \ref{assn:Main}  is equivalent to
\begin{assn} \label{assn:MainM}
For any point $p\in \sP$ and $t\in T^+$ above $p$, 
the eigensheaf $M$ satisfies  one of the following
conditions: 
\begin{itemize}
\item $M^\sP$ has only positive weights along the fiber $\pi^*(p)$ 
and $M_{\pi^*(p)}=N_{\pi^*(p)}$, or 
\item $0$ is a weight for $M^\sP$ and the support of the $0$-weight space is 
the point $t$.
\end{itemize}
\end{assn}

\subsection{Results}\label{ssec:results}

The paper is organized along the following lines: 
in Section \ref{nahmsect}, we give an overview of the conditions and results of \cite{Sz} which 
are most often referred to in the present paper. 

In Section \ref{sec:Basic}, we recall the notions which will be used throughout 
the paper: pure sheaves of dimension $1$, parabolic sheaves, parabolic Euler-characteristic, 
degree and stability of parabolic sheaves; and prove some of their properties. 

In Section \ref{sec:Blowup}, we define an {\em iterated} version of 
blow-up maps for non-reduced 
zero-dimensional subschemes. This will be essential for the 
generalization of Nahm transform 
to Higgs bundles with higher-order poles. 

In Section \ref{sec:ProperTransform}, analogous to the proper transform 
of a divisor with respect to  a blow-up, 
we introduce the {\em proper transform of a coherent sheaf} 
with respect to a blow-up of a closed point. We study properties of the proper transform for $1$-dimensional pure 
sheaves on surfaces. For such sheaves, proper transform is related to Hecke transforms of locally free sheaves. 
In particular, for such sheaves, proper transform is a quasi-inverse of the direct image 
(Lemma \ref{lem:qi}), and it preserves the Euler-characteristic (Lemma \ref{lem:EPpres}). 
We also give a parabolic version of the proper transform, and prove that it preserves the 
parabolic Euler-characteristic (Subsection \ref{subsec:ParPropTr}). 

In Section \ref{sec:AddDelete}, we define two operations to modify the divisor of 
parabolic sheaves: {\em Deletion} along $\sE$ removes an effective subdivisor $\sE$
of the parabolic divisor, whereas {\em addition} along $\sE$ appends an effective divisor $\sE$ to the parabolic divisor. 
Under an assumption (is equivalent to \ref{assn:Main}), these operations 
are inverse to each other. Moreover, they preserve the parabolic Euler-characteristic (Proposition \ref{prop:adddelinv}).

In Section \ref{sec:SpecTr}, we introduce what we call the {\em spectral triples}, consisting of
a smooth surface, an effective divisor on it, and a rank-one torsion-free 
sheaf on the divisor satisfying some properties. Notice that the operation $\pi^H$ 
associates to a given Higgs bundle a spectral triple $(Z^\sP,\supp(M^\sP),M^\sP)$. 
We shall call this spectral triple the {\em standard spectral triple} associated to 
the Higgs bundle. On the other hand, there exists another way of defining a spectral 
triple $(Z^0,\supp(M^0),M^0)$ where the surface is $Z^0=\CP \times \CPt$: we call this the 
{\em naive spectral triple}. The surfaces $Z^\sP$ and $\CP \times \CPt$ are related
by a series of elementary transformations. Let $Z$ be the resolution of 
indeterminacies of $Z^0 - \to Z^\sP$. Then, we show that the proper transforms of  $M^\sP$ and $M^0$ agree on $Z$ 
(Proposition \ref{prop:Relation}). 

In Section \ref{trsect}, we construct Nahm transform of parabolic Higgs bundles on the projective line
as a composition of the operations introduced up to this point. The starting point is the diagram 
\begin{equation} \label{diag:Nahm}
\xymatrix{
 & Z^{\abs} \ar[dl]_{\qmap_\sP} \ar[dr]^{\widehat 
 \qmap_{\widehat \sP}} \ar[d] & \\
 Z^\sP \ar[d] & \PP^1 \times \widehat \PP^1 \ar[dl] \ar[dr] 
 & \widehat Z^{\widehat \sP} \ar[d] \\
 \PP^1 & & \widehat \PP^1
}
\end{equation}
(see (\ref{absdiag})). Here the maps $\qmap_\sP$ and $\widehat{\qmap}_{\widehat \sP}$ are blow-up 
maps, and $Z^{\abs}$ is called the {\em intermediate spectral surface}. 
Starting from a $1$-dimensional parabolic sheaf $\M^\sP_\bullet$ on $Z^\sP$, 
Nahm transform produces a $1$-dimensional parabolic sheaf $\widehat\M^{\widehat\sP}_\bullet$ on 
$\widehat Z^{\widehat \sP}$ by the formula:
$$\M^\sP_\bullet \mapsto (-1)^*_{\widehat Z^{\widehat \sP}}\widehat\M^{\widehat\sP}_\bullet = 
(-1)^*_{\widehat Z^{\widehat \sP}}({\widehat \qmap_{\widehat \sP}})_*\add{\widehat \sE^+}\del{\sE^+}  \wtd{(\qmap_\sP)}
(\M^\sP_\bullet).$$
From right to left, this formula reads as 
a proper transform with respect to $\qmap_\sP$, deletion along a divisor $\sE^+$,
addition along a divisor $\widehat \sE^+$, push-forward with 
$\widehat \qmap_{\widehat \sP}$, and pull-back with respect to the fiberwise $(-1)$ multiplication. 
Here, $\sE^+$ and $\widehat \sE^+$ are suitably chosen divisors related 
to the birational 
morphisms $\qmap_\sP$ and ${\widehat \qmap_{\widehat \sP}}$ respectively.
Theorem \ref{trthm} shows that our construction generalizes that of \cite{Sz}.

In Section \ref{sec:Ex}, we describe two examples in which we use our method 
to compute  the transformed Higgs bundle explicitly. 
These examples are beyond the scope of \cite{Sz}. 
The first example features  a Higgs field with a nilpotent residue, whereas  
the second one  a higher-order pole. 

Section \ref{sec:qi} provides a geometric proof of 
the fact that the transformation 
is involutive up to a sign. 

In Section \ref{sec:moduli}, we study the map induced by Nahm transform on the 
moduli spaces of stable Higgs bundles of degree $0$ with prescribed singularity behaviour. 
First, we compute the dimension of these moduli spaces (Lemma \ref{lem:moddim}). 
Then, we show that Nahm transform preserves the parabolic degree, and for Higgs bundles of 
degree $0$, it preserves stability (Lemma \ref{stablem}). Finally, in Corollary \ref{isomcor} we prove that Nahm 
transformation induces a hyper-K\"ahler isometry between the corresponding moduli spaces.

\section{An Overview of analytic Nahm transform}\label{nahmsect}
In this section, we give a summary of the results of \cite{Sz} relevant for the present paper. 

Let $\sP=\{p_{1}, \ldots , p_{n} \}$ be a finite set in $\CP$ composed of distinct points at finite distance, 
$\E$ be a rank $r$ holomoprhic vector bundle on $\CP$ and 
$$
    \theta : \E \lra \E \otimes \O_{\CP}(\sP)
$$
be a holomorphic map (called the \emph{Higgs field}), where $\O_{\CP}(\sP)$ is the sheaf of meromorphic functions with 
at most simple poles in the points of $\sP$ and no other poles. We assume that in any $p_{j} \in \sP$ the Higgs
field has semi-simple residue: in the standard holomorphic coordinate $z$ of $\C$ and in a convenient 
holomorphic trivialization $\{ e^{j}_{1}, \ldots , e^{j}_{r} \}$ of $\E$ in a neighborhood of $p_{j}$ it can 
be written
\begin{equation}\label{thetaj}
   \theta =B_{j}\frac{1}{z-p_{j}} + O(1),
\end{equation}
where $O(1)$ stands for holomorphic terms and 
\begin{equation}\label{bj}
    B_{j}=\begin{pmatrix}
                  0 & & & & & \\
                  & \ddots & & & & \\
                  & & 0 & & & \\
                  & & & \lambda^{j}_{r_{j}+1} & & \\
                  & & & & \ddots & \\ 
                  & & & & & \lambda^{j}_{r}
           \end{pmatrix},
\end{equation}
is a diagonal matrix (the \emph{residue} of $\theta$ in $p_{j}$) with all the $\lambda^{j}_{k}$ for $r_{j}<k\leq r$
non-vanishing and distinct. We suppose furthermore that an  compatible \emph{parabolic structure} 
is given in $p_{j}$: this simply means the data of real numbers 
$0=\alpha^{j}_{1}=\ldots =\alpha^{j}_{r_{j}}<\alpha^{j}_{r_{j}+1}\leq \ldots \leq \alpha^{j}_{r}<1$ called 
\emph{parabolic weights}. Here, the condition $0=\alpha^{j}_{1}=\ldots =\alpha^{j}_{r_{j}}$ is an extra condition 
of admissibility. 
For any $\alpha \in [0,1[$ we then define the space $\F_{\alpha}\E_{p_{j}}$ to be the subspace of the fiber 
$\E_{p_{j}}$ of $\E$ in $p_{j}$ spanned by the $e^{j}_{k}(p_{j})$ such that $\alpha^{j}_{k}\geq \alpha$; this gives a 
finite filtration 
\begin{equation}\label{parj}
     \{ 0\} =\F_{1}\E_{p_{j}} \subset \F_{\alpha^{j}_{r}}\E_{p_{j}}\subset \cdots \subset \F_{\alpha^{j}_{r_{j}+1}}\E_{p_{j}}\subset \F_{0}\E_{p_{j}}=\E_{p_{j}}
\end{equation}
of the fiber. An alternative way to define a parabolic structure at a singularity on a Higgs bundle is to say 
that the holomorphic vector bundle has a parabolic structure (i.e. a filtration as above and parabolic weights) 
in the singularity, and the Higgs field is compatible with the parabolic structure in the sense 
that its residue endomorphism in the puncture is the sum of its induced endomorphisms on the graded vector
spaces of the parabolic filtration. 

At infinity, we suppose that $\theta$ is holomorphic, such that its Taylor series written in the local coordinate 
$z^{-1}$ and some holomorphic trivialization of $\E$ near infinity 
\begin{equation}\label{thetainf}
    \theta = \frac{1}{2}A+B_{\infty}\frac{1}{z}+ O\left(\frac{1}{z^{2}} \right)
\end{equation}
satisfy that the constant term $A$ is diagonal with eigenvalues $\xi_{1}, \ldots ,\xi_{\nt}$ of multiplicity 
possibly higher than one
\begin{eqnarray}\label{a}
     A =  \begin{pmatrix}
                \xi_{1} & & & & & & & \\
                & \ddots & & & & & & \\
                & & \xi_{1}  & & & & & \\
                & & & \ddots & & & & \\
                & & & & \ddots & & & \\
                & & & & & \xi_{\nt} & & \\
                & & & & & & \ddots & \\
                & & & & & & & \xi_{\nt}
                        \end{pmatrix}, 
\end{eqnarray}
and the first-order term $B_{\infty}$ is also diagonal (in the same trivialization) 
\begin{equation}\label{binf}  
     B_{\infty} =  \begin{pmatrix}
              \lambda^{\infty}_{1} & & & & & & & \\
               & \ddots & & & & & & \\
                & & \lambda^{\infty}_{a_{1}} & & & & & \\
                & & & \ddots & & & & \\
                & & & & \ddots & & & \\
                & & & & & \lambda^{\infty}_{1+a_{\nt}} & & \\
                & & & & & & \ddots & \\
                & & & & & & & \lambda^{\infty}_{r}
                        \end{pmatrix}. 
\end{equation}
Here the eigenvalues $\{ \lambda^{\infty}_{1+a_{l}},\ldots , \lambda^{\infty}_{a_{l+1}} \}$ correspond to the basis vectors spanning the 
$\xi_{l}$-eigenspace of $A$ (where we have put $a_{0}=0$). We make the assumption that for a fixed $1\leq l\leq \nt$, 
none of these eigenvalues $\{ \lambda^{\infty}_{1+a_{l}},\ldots , \lambda^{\infty}_{a_{l+1}} \}$ vanishes and they are all distinct. 
We suppose that a parabolic structure is given in this singularity as well: that is, we are given parabolic
weights $0<\alpha^{\infty}_{k}<1$ for $k=1,\ldots ,r$, arranged in such a way that inside one block $a_{l}<k\leq a_{l+1}$ 
they form an increasing sequence, and a corresponding filtration $\F_{\alpha}\E_{\infty }$ for $\alpha \in [0,1]$ of the fiber 
of $\E$ over infinity, spanned by the basis elements having parabolic weight $\geq \alpha$.  
A sheaf $\E$ with a parabolic structure will often be denoted by $\E_{\bullet}$.

\begin{rk}
Multiplication by the globally defined meromorphic $1$-form $\d z$ associates to $\theta$ an endomorphism-valued 
meromorphic $1$-form that has a double pole at infinity because the form $\d z$ does. Therefore, we will 
often call $A/2$ the \emph{second-order term} and and $B_{\infty}$ the \emph{residue} of $\theta$ at infinity.  
Also, a holomorphic vector bundle with a parabolic structure in $\sP\cup \{\infty \}$ will be called \emph{parabolic 
vector bundle}; if moreover a compatible Higgs field is given, then we will call it a \emph{parabolic 
Higgs bundle}. 
\end{rk}

Denote by $\deg(\F)$ the usual algebraic geometric degree of a holomorphic vector bundle $\F$. It is 
clear that a {\em subbundle} $\F$ (or quotient bundle $\Q$) of a parabolic vector bundle $\E$ also admits 
an induced parabolic structure by intersecting with $\F$ the terms of the filtration of $\E$, 
and assigning the biggest of the weights to all filtered terms that become isomorphic after taking 
intersections with $\F$. 

\begin{defn}
The \emph{parabolic degree} of $\E_{\bullet}$ is the real number 
$\para\deg(\E_{\bullet} )=\deg(\E) + \sum_{j\in \{1,\ldots ,n,\infty \}}\sum_{k=1}^{r}\alpha^{j}_{k}$. The \emph{parabolic slope}
of $\E_{\bullet}$ is the real number $\para \mu(\E_{\bullet} )= \para\deg(\E_{\bullet})/rk(\E)$. Finally, $(\E_{\bullet},\theta )$ is 
said to be \emph{parabolically stable} if for any subbundle $\F_{\bullet}$ invariant with respect to $\theta$ with 
its induced parabolic structure, the inequality $\para \mu(\F_{\bullet}) < \para \mu(\E_{\bullet})$ holds. 
\end{defn}

Suppose in all what follows that $(\E,\theta)$ is not the trivial line bundle $\O_{\CP}$ together with a 
constant multiplication map. 
Denote by $\Ct$ the \emph{dual line} of $\C$ (another copy of $\C$), and by $\CPt$ the \emph{dual sphere}, 
the compactification of $\Ct$ by the point $\ti$. By \cite{Sz}, 
the Nahm transform of a stable parabolic Higgs bundle $(\E_{\bullet} , \theta )$ of parabolic degree $0$ is then a parabolic Higgs bundle 
$(\Et_{\bullet}, \tt)$ on $\CPt$, with regular singularities (i.e. $\tt \d \xi$ having simple poles) in the set 
$\Pt = \{ \xi_{1}, \ldots ,\xi_{\nt} \}$ and an irregular singularity (i.e. $\tt \d \xi$ having a double pole, therefore 
$\tt$ being holomorphic) at infinity. Also, the transform of a Hermitian-Einstein metric on $(\E_{\bullet} , \theta )$ 
is a Hermitian-Einstein metric on $(\Et_{\bullet}, \tt)$; in particular, this
latter is poly-stable. We sketch the idea 
of the construction of the transform. First, introduce a twist of the Higgs field: for any $\xi \in \Ct$ set 
\begin{equation*}
    \theta_{\xi} = \theta - \frac{\xi}{2} \mbox{Id}_\E
\end{equation*}
where $\mbox{Id}_\E$ is the identity bundle endomorphism of $\E$. 
Consider now the \emph{open spectral curve} $\Sigma^{\oo}$ in $(\CC \setminus \sP) \times (\Ct \setminus \Pt)$ defined by 
$$
    \Sigma^{\oo}=\{ (z,\xi ) \: | \: \det(\theta_{\xi})(z)=0 \}. 
$$
In other words, denoting by $\pi^{\oo}$ (respectively $\pit^{\oo}$) the projection on $\CC \setminus \sP$ 
(respectively $\Ct \setminus \Pt$) in the product $(\CC \setminus \sP) \times (\Ct \setminus \Pt)$, 
this curve is the support of the cokernel sheaf $M^{\oo}$ of the map 
$$
    \theta_{\xi} : (\pi^{\oo})^{\ast}\E \lra (\pi^{\oo})^{\ast}\E.
$$
Propositions 4.7, 4.15 and 4.24 together with Lemma 4.30 of \cite{Sz} give an interpretation of
the transform on an open set. 
\begin{thm}\label{openint} The Nahm transformed Higgs bundle restricted to $\Ct \setminus \Pt$ can be obtained as follows: 
\begin{itemize}
\item{the holomorphic bundle $\Et$ is the pushdown $\pit^{\oo}_{\ast} M^{\oo}$ endowed with its induced 
holomorphic structure; we denote its rank by $\rt$} 
\item{on the open set of $\xi \in \Ct \setminus \Pt$ over which the fiber of $\Sigma^{\oo}$ consists of distinct points 
$\{z_{1}(\xi), \ldots , z_{\rt}(\xi)\}$ of multiplicity $1$, the transformed Higgs field $\tt$ acts on the subspace 
$\coker(\theta_{\xi}(z_{k}(\xi))) \subset \Et|_{\xi}$ as multiplication by $-z_{k}(\xi)/2$; this then admits a unique continuation
into points where the fiber has multiple points.} 
\end{itemize}
\end{thm}
This description then gives an understanding of the behaviour of the Higgs field near a point of 
$\Pt$ and near $\ti$: we only have to understand the behaviour of the open spectral curve near 
these points. Because of the special form of $\theta$ in the singularities, we deduce that 
the eigenvalues of the transformed Higgs field have indeed simple poles in the points of 
$\Pt$, and are bounded near $\ti$. In different terms, this defines a natural compactification 
$\Sigma^{0} \subset \CP \times \CPt$ of $\Sigma^{\oo}$. Moreover, we gain precise information about its asymptotic 
expansions near these points: namely, near a point $\xi_{l} \in \Pt$ the residue of the transformed Higgs 
field in a convenient trivialization of the transformed bundle is equal to 
\begin{eqnarray*}
           -\begin{pmatrix}
                  0 & & & & & \\
                  & \ddots & & & & \\
                  & & 0 & & & \\
                  & & & \lambda^{\infty}_{1+a_{l}} & & \\
                  & & & & \ddots & \\ 
                  & & & & & \lambda^{\infty}_{a_{l+1}}
           \end{pmatrix},
\end{eqnarray*}
in other words it is the direct sum of the opposite of the residue of the original Higgs field at infinity 
restricted to the $\xi_{l}$-eigenspace of the leading term and a $0$-matrix 
(Theorem 4.32 of \cite{Sz}); whereas its leading term at infinity in a convenient trivialization is 
\begin{eqnarray*}
              -\frac{1}{2} \begin{pmatrix}
                p_{1} & & & & & & & \\
                & \ddots & & & & & & \\
                & & p_{1}  & & & & & \\
                & & & \ddots & & & & \\
                & & & & \ddots & & & \\
                & & & & & p_{n} & & \\
                & & & & & & \ddots & \\
                & & & & & & & p_{n}
                        \end{pmatrix},
\end{eqnarray*}
each $p_{j}$ appearing with multiplicity $rk(res(\theta,p_{j}))=r-r_{j}$, and the corresponding first-order term in
the same trivialization is then 
\begin{eqnarray*}
             -\begin{pmatrix}
                \lambda^{1}_{r_{1}+1} & & & & & & & \\
                & \ddots & & & & & & \\
                & & \lambda^{1}_{r} & & & & & \\
                & & & \ddots & & & & \\
                & & & & \ddots & & & \\
                & & & & & \lambda^{n}_{r_{n}+1} & & \\
                & & & & & & \ddots & \\
                & & & & & & & \lambda^{n}_{r}
                        \end{pmatrix}
\end{eqnarray*}
(Theorem 4.33 of \cite{Sz}). In particular, we deduce the formula 
\begin{equation}\label{trrank}
    \rt = \sum_{j=1}^{n} rk(res(\theta,p_{j})).
\end{equation}
Therefore, we see an intricate interplay between singularity behaviour at the regular 
singularities and the one at the irregular singularity. 

Afterwards, we use the extensions of $\E$ over the singularities to define an extension $M^{0}$ of 
$M^{\oo}$ to the compactified spectral curve $\Sigma^{0}$. These in turn induce an extension
$\Et_{\bullet}^{ind}$ of $\Et$ into a holomorphic bundle endowed with a parabolic structure in each point 
of $\Pt \cup \{ \ti \}$, which we call the \emph{induced extension} (c.f. \cite{Sz}, Section 4.4). 
By definition, a local holomorphic section of this extension has 
a $D''_{\xi}$-harmonic representative obtained from a local section of the cokernel sheaf $M^{0}$ 
multiplied with a bump-function of constant height concentrated near the spectral points of $\xi$,
such that the diameter of their support converges to $0$ up to first order near the points $\xi_{l}$ and 
to $\infty$ near $\ti$ also up to first order. Next, we compute the parabolic weights of these extensions 
with respect to the transformed Hermitian-Einstein metric: for a point $\xi_{l} \in \Pt$ the non-zero weights 
are equal to $\alpha^{\infty}_{1+a_{l}}-1, \ldots ,\alpha^{\infty}_{a_{l+1}}-1$; whereas the weights at $\ti$ are equal to 
$\alpha^{1}_{r_{1}+1}-1, \ldots ,\alpha^{1}_{r}-1, \ldots , \alpha^{n}_{r_{n}+1}-1, \ldots ,\alpha^{n}_{r}-1$ (c.f. \cite{Sz}, Section 4.6). 
In particular, all the non-zero weights violate the requirement that they be between $0$ and $1$: 
they are actually shifted by $-1$. 
Therefore, in order to get a genuine parabolic Higgs bundle on $\CPt$, we have to change the induced 
extension at $\ti$ by a factor of $\xi^{-1}$, and the extension of the basis vectors corresponding to non-zero 
eigenvalues of the Higgs field at the logarithmic singularities $\xi_{l}$ by factors 
of $(\xi-\xi_{l})$. The result we obtain this way is called the \emph{transformed extension}, and denoted
$\Et^{tr}$ (c.f. \cite{Sz}, Section 4.7). 
Finally, an application of Grothendieck-Riemann-Roch theorem yields the degree of the transformed 
holomorphic bundle on $\CPt$ with respect to the transformed extensions:
\begin{equation}\label{degformula}
    \para \deg(\Et_{\bullet}^{tr} ) = \para \deg(\E_{\bullet} ).
\end{equation}

\section{Basic Material}\label{sec:Basic}

Fix  a projective scheme $X$ over a field $\kk$ with an ample invertible
sheaf $\O_X(1)$.
For a given coherent $\O_X$-module $\E$, 
the {\em support} of $\E$ is the closed set $\supp(\E)=\{ x\in X | \, \E_x\neq 0\}$.
Its dimension is called the {\em dimension} of the sheaf $\E$ and is denoted by $\dim(\E)$. 

\begin{defn}
A coherent sheaf $\E$ of $\O_X$-modules on $X$ is {\em pure of dimension $d$} if  $\dim(\F)=d$ for any nontrivial 
subsheaf $\F$ of
$\E$.
\end{defn}

\begin{defn}
A subsheaf $\F$ of a pure $d$-dimensional sheaf $\E$ is {\em saturated} 
if $\E/\F$ is either $0$ or pure of dimension $d$.
\end{defn}

Equivalently, $\E$ is pure if and only if the associated points of $\E$ are all of the same dimension.
\begin{defn}
For a given coherent sheaf $\E$ on $X$,
the {\em Euler characteristic} is defined to be 
$\chi(\E)=\sum_{i=0}^{\dim(X)} (-1)^i \dim_k H^i(X, \E)$. The {\em Hilbert polynomial}
$P(\E)$ of $\E$ is defined by $P(\E,m):=\chi(E\otimes \O_X(m))$.
\end{defn}

Let $\sD$ be a Cartier divisor on $X$. Suppose $\E$ is a pure sheaf of dimension
$d$ and $\dim(\sD\cap\supp\E)<\dim\supp\E$, then $\E\otimes\O_X(-\sD)\ra \E$
is injective.

\begin{defn} \label{def:parsh}
Let $\E$ and $\sD$ be as above.
A triple $(\E,F_\bullet\E,\alpha_\bullet)$ is called a {\em parabolic 
sheaf} on $X$ with parabolic divisor $\sD$ and weights $\alpha_\bullet$ if 
$F_\bullet\E$ is a filtration of  $\E$ by coherent subsheaves $F_i\E$ so that
$
\E(-\sD)=F_{l+1}\E\subset
F_l\E \subset \cdots \subset F_2\E \subset F_1\E=\E
$
and $\alpha_\bullet$ is  a sequence of real numbers $0\leq \alpha_1 \leq \alpha_2 \leq \cdots \leq \alpha_l<1$.
Set $\gr{F}{i}\E:=F_i\E/F_{i+1}\E$.
\end{defn}
One can view  $\gr{F}{i}\E$ as coherent sheaves on $\sD\cap \supp(\E)$.

\begin{defn}
Let $\sD_1$ be an irreducible component of $\sD$ such that for any other irreducible component $\sD'$ one has 
$\sD_1\cap \sD' \cap \supp \E= \emptyset$. We say that the parabolic structure is trivial on $\sD_1$ if 
$\gr{F}{i}\E|_{\sD_1}=0$ for all $i>1$ and $\alpha_1=0$.
\end{defn}

\begin{defn} Let $\E$ and $\sD$ be as before.
The pair $(\E,\E_\bullet)$ is an {\em $\RR$-parabolic sheaf} on $X$ if  $\E_\bullet=\{\E_\alpha\}$ 
 is a collection of
coherent sheaves
parametrized by $\alpha\in\RR$ satisfying the following
properties:
\begin{enumerate}
\item $\E_0=\E$,
\item For all $\alpha<\beta$, $\E_\beta$ is a coherent subsheaf of $\E_\alpha$,
\item For all $\alpha$ and small $\eps>0$, $\E_{\alpha-\eps}=\E_\alpha$,
\item For all $\alpha$, $\E_{\alpha+1}=\E_\alpha(-\sD)$.
\end{enumerate}
Set $\gr{F}{\alpha}:=\E_\alpha/\E_{\alpha+\eps}$ for small $\eps>0$.
\end{defn}

Parabolic sheaves and $\RR$-parabolic sheaves are equivalent: To see this,
set $\alpha_{l+1}=1$ and $\alpha_0:=\alpha_l-1$. For any real number $\alpha$,
let $i$ be the unique integer so that $\alpha_{i-1}<\alpha-\floor{\alpha} \leq \alpha_i$
where $\floor{\alpha}$ is the largest integer with $\floor{\alpha}\leq \alpha$.
Set $\E_\alpha:=F_i\E(-\floor{\alpha}\sD)$. Conversely, given an $\RR$-parabolic
sheaf $\E_\bullet$, inductively choose $0\leq \alpha_i<1$ for $i=1,\ldots,l$ 
so that 
$\E_{\alpha_i}$
properly contains $\E_\beta$ for any $\beta>\alpha_i$.
Set $F_i\E:=\E_{\alpha_i}$ and $F_{l+1}=\E(-\sD)$. The resulting triple 
$(\E, F_\bullet, \alpha_\bullet)$ is
a parabolic sheaf on $X$.  
Therefore, when the ($\RR$-)parabolic structure of a sheaf $\E$ is clear from
the context, we will write $\E_\bullet$ for the pair $(\E,\E_\bullet)$.

\begin{defn}
Given two parabolic sheaves $\E'_\bullet$ and $\E_\bullet$, an $\O_X$-module
homomorphism $\varphi:\E'\ra \E$ is a {\em parabolic homomorphism} 
if $\varphi(\E'_\alpha)\subset \E_\alpha$ 
for all real numbers $\alpha$.
For $\E$ a parabolic sheaf and $\E$ a subsheaf endowed with a parabolic structure, 
we say that $\E'_\bullet$ is a {\em parabolic subsheaf} if 
the inclusion is a parabolic homomorphism.
\end{defn}

\begin{defn}
Given a saturated subsheaf $\E'$ of a parabolic sheaf $\E_\bullet$, the {\em induced parabolic structure}
$^{ind}\E'_\bullet$ on
$\E'$ is  defined as $^{ind}\E'_\alpha:=\E'\cap \E_\alpha$ for all $\alpha\in \RR$. 
\end{defn}
\begin{rk}\label{rk:indparstr}
The induced parabolic structure
for a given saturated subsheaf $\E'$ is the largest among the parabolic structures $\E'_\bullet$ which 
make $\E'$ into a parabolic subsheaf of $\E_\bullet$. As a consequence,  
it suffices to consider the saturated subsheaves  of a parabolic sheaf
$\E_\bullet$ 
with their induced 
parabolic structures
to 
measure the stability of $\E_\bullet$. 
\end{rk}

\begin{defn}
The {\em parabolic Euler characteristic}
of a 
parabolic sheaf $\E_\bullet$ is defined as follows:
\begin{eqnarray}
\para\chi(\E_\bullet) &:= &
\chi(\E(-\sD))+\sum_{i=1}^l \alpha_i\chi(\gr{F}{i}\E) 
\end{eqnarray}

If $X$ is a curve, then the {\em parabolic degree} of $\E_\bullet$ is defined as: 
\begin{eqnarray}
\para\deg(\E_\bullet)&:=&\deg(\E )+ \sum_{i=1}^l \alpha_i \dim(\gr{F}{i}\E)
\end{eqnarray}

\end{defn}

One can check that $\para\chi(\E_\bullet)=\int_0^1\chi(\E_\alpha)\,d\alpha$ (see \cite{yo}).


\begin{prop}
The parabolic Euler characteristic is additive: Given any short exact sequence 
$$
   0\longrightarrow \E'_\bullet \longrightarrow \E_\bullet \longrightarrow \E''_\bullet \longrightarrow 0
$$
of parabolic sheaves with the same parabolic divisor $\sD$, the identity $$\para\chi(\E_\bullet)=\para\chi(\E'_\bullet)+\para\chi(\E''_\bullet)$$ holds. 
\end{prop}

\begin{proof}
Recall that a sequence of parabolic sheaves is said to be exact if for all $\alpha \in \RR$ the induced sequence on the $\alpha$-filtered 
terms is exact. Taking $\alpha =-1$, we see that 
$$
   0\longrightarrow \E'(-\sD) \longrightarrow \E(-\sD) \longrightarrow \E''(-\sD)\longrightarrow 0
$$
is exact. By additivity of the usual Euler characteristic, $\chi(\E(-\sD))=\chi(\E'(-\sD))+\chi(\E''(-\sD))$. 
On the other hand, the snake lemma implies that for any $\alpha \in \RR$ the induced sequence on the $\alpha$-graded pieces 
$$
   0\longrightarrow \gr{F}{\alpha} \E' \longrightarrow \gr{F}{\alpha}\E \longrightarrow \gr{F}{\alpha}\E''\longrightarrow 0
$$
is also exact. The statement follows by applying additivity of $\chi$ to these sequences. 
\end{proof}

\begin{defn}
Given a parabolic sheaf  $\E_\bullet$ and $L$ a line bundle, 
define a parabolic structure on $\E\otimes L$ by 
setting $(\E \otimes L)_{\alpha}:=\E_{\alpha} \otimes L$ for all $\alpha\in \RR$. 
\end{defn}

When $L=\O(\sD)$,
this definition of the parabolic structure of $\E(\sD)$ is to some extent unnatural from the analytic point of view:
indeed, the behaviour of a fixed harmonic metric on local sections of $\E(\sD)$ is not the same as on local sections 
of $\E$, but changes by a factor of $|z|^{-2}$ where $z$ is a local coordinate centered at $\sD$. 
However, we have two reasons to choose this convention: 
first, the proposition above which says that the parabolic  Euler characteristic is additive if the weights of all 
parabolic sheaves are in the same interval; second, for any parabolic Higgs sheaf $(\E,\theta)$ (see the definition below)
with divisor $\sD$, this definition makes $\theta :\E \rightarrow \E(\sD)$ a morphism of parabolic sheaves.

On a smooth projective curve $X$, parabolic Euler characteristic and parabolic
degree are related as follows:
\begin{prop} \label{prop:parRR}
If $\E$ is a parabolic sheaf on a smooth projective curve $X$, then 
\begin{align*}
   \para\deg(\E(\sD))&=\para\deg(\E) + r\deg{\sD},\\
    \para\chi(\E(\sD))&=\para\deg(\E) + r\chi(\O_X).
  \end{align*}
\end{prop}

\begin{proof}
The first formula follows by definition, because the jumps $\alpha_i$ of the parabolic structures 
of $\E$ and $\E(\sD)$ are the same, and the graded pieces of the filtration corresponding to each  $\alpha_i$ 
are isomorphic. 

The second follows from Riemann-Roch and the isomorphism of the graded pieces: 
\begin{align*}
   \para\chi(\E(\sD))&=\chi(\E)+\sum_{i=1}^l \alpha_i\chi(\gr{F}{i}\E(\sD))\\
		    &=\deg(\E)+r\chi(\O_X)+\sum_{i=1}^l \alpha_i\chi(\gr{F}{i}\E)\\
		    &=\deg(\E)+r\chi(\O_X)+\sum_{i=1}^l \alpha_i\dim(\gr{F}{i}\E)
\end{align*}
because the $\gr{F}{i}\E$ are supported on the $0$-dimensional subscheme $\sD$. 
\end{proof}

\begin{defn} \label{defn:Higgssheaf}
A {\em Higgs sheaf}
$(\E,\theta)$ consists of 
a coherent sheaf $\E$ on $X$ together with a $\O_X$-module homomorphism $\theta:\E \ra \E(\sD)$. 
The resulting $\O(\sD)$-valued endomorphism $\theta$ is called a {\em Higgs field}.
A {\em parabolic Higgs sheaf} $(\E_\bullet,\theta)$ with divisor $\sD$
consists of 
a parabolic sheaf $\E_\bullet$ on $X$ with divisor $\sD$ and a parabolic homomorphism $\theta:\E_\bullet \ra \E_\bullet(\sD)$. 
A {\em homomorphism of Higgs sheaves} $\psi: (\E^1,\theta^1) \ra (\E^2,\theta^2) $
is homomorphism of sheaves $\psi: \E^1 \ra \E^2$ commuting with the Higgs fields:
$(\psi\otimes 1)\circ\theta^1=\theta^2\circ\psi$. A {\em homomorphism of
parabolic Higgs sheaves} is a homomorphism of Higgs sheaves respecting the
parabolic structure.

\end{defn}
A (parabolic) Higgs subsheaf of $\E$ is defined  in the obvious way: it is a (parabolic) subsheaf preserved by the Higgs field.

\begin{rk}\label{rk:divisors}
Starting from Section \ref{trsect}, we will consider Higgs sheaves on $X=\PP^1$ 
with polar divisor $\sP$ and parabolic divisor $\sD =\sP +\infty$. In terms of Definition 
\ref{defn:Higgssheaf}, these objects are defined as Higgs sheaves with polar divisor $\sD$ 
with an apparent singularity at $\infty$. In other words, the Higgs field, as a rational 
section with values in $\sP$, extends regularly at $\infty$. We take $Z^\sP$ as standard spectral 
surface for a Higgs sheaf with polar divisor $\sP$ and parabolic divisor $\sD$. It would also 
be possible to work with the surface $Z^\sD$ -- these two surfaces are related by an elementary 
transformation over the infinity. However, we work with $Z^\sP$ because the poles of the Higgs field
are already contained in $\sP$. 
\end{rk}

\begin{defn}
A parabolic (Higgs) sheaf $\E_\bullet$ 
is said to be {\em semistable} if for any given
proper parabolic (Higgs) 
sheaf $\F_\bullet \subset \E_\bullet$, $\para p(\F_\bullet,m) \leq \para p(\E_\bullet,m)$ for large $m$. The (Higgs) sheaf $\E_\bullet$ is said to
be {\em stable} if for all proper parabolic (Higgs) subsheaves 
$\F_\bullet \subset \E_\bullet$, $\para p(\F_\bullet,m) <\para p(\E_\bullet,m)$ for large $m$.
\end{defn}

\subsubsection{} \label{sec:ParabolicHiggsM} The standard construction described in 
Section \ref{ssec:HiggsM} adapts to the parabolic case as well.
A parabolic Higgs sheaf $\theta:\E_\bullet\ra\E_\bullet(\sP)$ with 
parabolic divisor $\sD$ determines a parabolic sheaf $M^\sP_\bullet$ on $Z^\sP$ with
parabolic divisor $\pi_\sP^*(\sD)$, with
$\pi_{\sP*}M^\sP_\alpha=\E(\sP)_\alpha$ for any $\alpha\in \RR$ and $\supp M^\sP_\bullet \cap
(y_\sP)=\emptyset$. Write $\pi^H(\E_\bullet,\theta)$ for $M^\sP_\bullet$ and
$\pi_H(M^\sP_\bullet)$ for $(\E_\bullet,\theta)$.

\subsubsection{Automorphism $(-1)$} \label{sec:(-1)}
If $M$ corresponds to the Higgs sheaf $\theta:\E\ra \E(\sP)$,
the pullback $(-1)_{Z^\sP}^*M$ 
corresponds to $-\theta:\E \ra \E(\sP)$.
We formalize this for parabolic Higgs sheaves as well as Higgs sheaves:

\begin{lem} \label{lem:(-1)}
\begin{align*}
\pi_H\, (-1)^*_{Z^\sP} \, \pi^H(\E,\theta) & =\pi^H(\E,-\theta),\\
\pi_H\, (-1)^*_{Z^\sP} \, \pi^H(\E_\bullet,\theta) & =\pi^H(\E_\bullet,-\theta).
\end{align*}
\end{lem}

\subsection{Commutative Algebra}
Let A be a local ring with maximal ideal $\mf$. For any $A$-module $M$, 
the depth of $M$ is defined as 
$$\depth(\M):=\min\{ i : \extgl{i}{A}{A/\mf}{M} \neq 0 \}$$
and the homological dimension $\dim(\M)$ is defined as the minimal 
length of a projective resolution of $\M$.

The Auslander-Buchsbaum formula relates the two invariants:
$$\dimh(\M)+\depth(\M)=\depth(A).$$
If $A$ is a regular local ring, then $\depth(A)=\dim(A)$.

\begin{cor} \label{cor:AusBuchs}
Let $A$ be a regular local ring of dimension $2$, $\M$ a torsion $A$-module with
$\dimh(\M)=1$. Then, any submodule $\M'$ of $\M$ is a torsion $A$-module with
$\dimh(\M')=1$.
\end{cor}
\begin{proof}
Any submodule $M'$ of $M$ is torsion, therefore not locally free and 
$\dimh(\M')\geq 1$. 
Since $\depth(\M)=1$, 
$$\extgl{0}{A}{A/\mf}{\M'}\subset \extgl{0}{A}{A/\mf}{\M}=0.$$
Consequently, $\depth(\M')\geq 1$. By Auslander-Buchsbaum equality,
$$\dimh(\M')=\depth(\M')=1.$$
\end{proof}

\begin{lem}
Assume $X$ is a smooth projective surface. For a coherent sheaf $M$ on $X$ ,
the following are equivalent:
\begin{enumerate}
\item $\M$ is pure of dimension $1$,
\item $\M$  is a torsion sheaf with $\dim(\M_x)=1$ for all $x\in X$,
\item $\dimh(\M_x)=\depth(\M_x)=1$ for all $x\in X$.
\end{enumerate}
Moreover, any subsheaf of a given pure sheaf $M$ of dimension $1$ is also pure
of dimension $1$. 
\end{lem}

\begin{proof}
Apply Cor \ref{cor:AusBuchs} and Prop. 1.1.10 \cite{HL}, which in this particular
case, states that a coherent sheaf $M$ of dimension $1$ is pure if and only if 
$\depth(\M_x)\geq 1$ for all $x\in X$.
\end{proof}

\section{Iterated Blow-Ups}\label{sec:Blowup}
A {\em sequence} of infinitesimally near points 
$(p_0,\ldots, p_n)$
on $X$ is defined recursively as follows: 
Let $X_0:=X$ and $p_0\in X_0$. 
By $\blow_j:X_j \lra X_{j-1}$,
denote the blow-up of $X_{j-1}$ at $p_{j-1}$, the exceptional divisor
$\blow_j^{-1}(p_{j-1})$ by $\sE_j$
and let $p_j$ be a point in $\sE_j$ 
for $j=1,\ldots, n$. By abuse of notation, denote the total transform of 
the exceptional divisor $\sE_j$ in $X_n$ still by $\sE_j$ 
for $j=1,\ldots, n$. For $1<j<n$, 
set
\begin{equation} \sC_j:=\sE_j-\sE_{j+1}.\end{equation}
The curve $\sC_j$ is a $(-2)$-curve on $X_n$.  
 
\begin{defn}
We call a zero dimensional closed subscheme $T$ of a smooth
(projective) surface $X$ {\em linear} if for each $p \in T_{red}$,
one can find $u,v\in \mf_{X,p}$ and positive integer $n$ so that 
$$ \mf_{X,p}=(u,v)\O_{X,p} \qquad \textnormal{ and } \qquad 
\J_{T,p}=(u^n,v)\O_{X,p}.$$
\end{defn}

For $p \in T_{red}$, the integer $n$ is uniquely determined and is equal 
to the length of $T$ at $p$, $\dim_k(\O_{X,p}/\J_{T,p})$. Denote this
integer by $n_p$. The total length $N$ of $T$ equals the sum $\sum n_p$. 

An {\em irreducible} linear subscheme $T$ of local length $n+1$ with closed point 
$p$ determines 
a sequence of infinitesimally near points $(p_0,\ldots,p_n)$ with $p_0=p$ as 
follows: Let 
\begin{itemize}
\item $p_0:=p$, 
\item $\sD_0:=(u^n-v)$, 
\item $\sD_j:=\wtd \blow_j \sD_{j-1}$ and
\item $p_j$ be the unique intersection point of $\sD_j$ with $\sE_j$ for 
$j=1,\ldots, n$.
\end{itemize}

The divisor $\sD_0$ is a smooth curve, thus so are all $\sD_j$ for $j>0$.  Because $mult_{p_0}\sD_0=1$,
it follows that $\sD_j\cdot \sE_j=1$, i.e. the intersection of $\sD_j$ and $\sE_j$ is unique  point, say $p_j$, for $j>0$.


We call the surface $X_n$ the iterated blow-up of $X$ at $T$ and
denote it by $\blow_T:\fbl{T}{X} \lra X$.


Enumerate the components of a linear subschemes as $T_1,\ldots, T_m$.
Then define the iterated blow-up of $X$ at $T$ to be

$$\fbl{T}X:=\fbl{T_1}X \times_X\cdots\times_X\fbl{T_m}X$$
and $\blow_T:\fbl{T}X \lra X$ be the corresponding morphism.

Clearly, $\fbl{T}X$ and $\blow_T$ do not depend on the enumeration chosen.
However, we need the enumeration for better record keeping: Denote
the closed point of $T$ corresponding to $T_i$ by $p_i$ and add the
subscript $i$ in front of previously written subscripts 
for the related data, thus making them 
$p_{ij}, \sE_{ij}, \sC_{ij}$
for appropriate values of $j$.

\subsection{Formulas for Exceptional Divisors}
Each leg of the following diagram 
is an iterated blow-up. To keep the notation simpler, assume that
$\sD=\n\cdot \textnormal{pt}$ for some $\n>0$ and replace $\sD$ with $\n$ in notation,
making $Z^\sD$ into $Z^\n$ etc.
$$
\xymatrix{
& \ar[ld]_{\pmap_0} Z  \ar[rd] ^{\pmap_n} & \\
Z^0  \ar@{-->}[rr] && Z^\n}
$$
In order to construct this diagram, one has to fix a global section $s$ of $\O(\sD)$
so that $\sD=\div{s}$. All such divisors differ by non-zero multiple. Let $u$ be 
a global section of $\O(\mathrm{pt})$, without loss
generality assume that $s=u^n$.

Given divisor in $Z^0$ or $Z^n$, denote its total transform in $Z$ by the same
letter.
Attach a superscript $+/-$ to divisors related 
to $\pmap_n$ and $\pmap_0$ respectively. 
Denote the fiber class by $\sF$ on any of the surfaces $Z^0$, $Z^n$ and $Z$.
Moreover, set
\begin{eqnarray}
\sC_0^\pm & := & \sF - \sE_1^\pm, \\
\sC_n^\pm & := & \sE_n^\pm.
\end{eqnarray}
Recall that
\begin{eqnarray}
\sC_j^\pm & = & \sE_j^\pm - \sE_{j+1}^\pm \quad\textnormal{for}\quad j=1,\ldots, n-1. 
\end{eqnarray}
Then,\begin{eqnarray}
\sF & = & \sE_j^+ + \sE_k^-  \quad\textnormal{for}\quad j+k=n+1,\\
\sF & = & \sum_{i=0}^n \sC_i^+ = \sum_{i=0}^n \sC_i^-.  
\end{eqnarray} 

Denoting the linear equivalence of divisors $\sim$, we see that
\begin{eqnarray}
\sX_0 & \sim & \sY_0, \\
\sX_n & \sim & \sY_n+n\cdot \sF. 
\end{eqnarray}

The formulae below relate various (exceptional) divisors. 
\begin{lem}
\begin{eqnarray}
\sX_0 & = & \sX_n -\sum_{j=1}^n \sE_j^+ \\
\sY_n & = & \sY_0 - \sum_{j=1}^n \sE_j^- \\
\sC_j^- & = & \sC_{n-j}^+ \quad\textnormal{for}\quad j=0,\ldots, n 
\end{eqnarray}

This table summarizes various relations:
\begin{equation}
\begin{array}{ccccccc}
\sE_n^- &=& \sC_n^- &=& \sC_0^+ &=& \sF-\sE_1^+ \\
\sE_{n-1}^- - \sE_n^- &=& \sC_{n-1}^- &=&  \sC_1^+ &=&  \sE_1^+-\sE_2^+ \\
\vdots & &\vdots& & \vdots& & \vdots \\
\sE_1^--\sE_2^- &=&  \sC_1^- &=&  \sC_{n-1}^+ &=&  \sE_{n-1}^+-\sE_n^+ \\
 \sF-\sE_1^- &=&  \sC_0^- &=&  \sC_n^+ &=&  \sE_n^+ 
\end{array}
\end{equation}
\end{lem}

We switch from the additive notation of divisors to multiplicative notation
of line bundles and sections: Let $\sF$ be the fiber above $\mathrm{pt}$,
i.e. it is cut out by
the equation $u=0$. Denote the section corresponding to divisor by the same letter
in small case,
i.e. $x_0$ is the section which cuts the divisor $\sX_0$. 
Set $\sC_i=\sC_i^+$. The equality 
$\sF=\sum \sC_i^+$, now becomes $u=\prod_{0}^n c_i$. 

Given a Higgs bundle $\theta: \E \lra \E(\sD)$, the eigensheaf
$\M^\sD$ on $Z^\sD$ equals $\coker(x_n-y_n\theta)$. We want to relate 
$\M^\sD$ to $\M^0$ on $Z^0$. The sheaf $M^0$ is defined by using $sx_0-y_0\theta$.
We  relate $sx_0-y_0\theta$ to $x_n-y_n\theta$.

\begin{eqnarray}
sx_0 - y_0 \theta & = & u^nx_0-y_0\theta \nonumber \\
& = & (\prod_1^{n}e_j^-)(x_n-y_n\theta) \nonumber \\
& = & (\prod_0^{n-1}c_i^{n-i})(x_n-y_n\theta) 
\end{eqnarray}

For $0\leq k \leq n$, 
\begin{eqnarray*}
x_n & = & (\prod_{1}^{n} e_j^+) x_0 \nonumber \\
 & = & (\prod_{1}^{n} c_j^{n+1-j}) x_0 \\
\gcd(x_n,u^k) & = & c_1^1c_2^2\cdots c_{k-1}^{k-1} c_k^k \cdots c_n^k \\
(\prod_{1}^{n}e_j^-)\gcd(x_n,u^k) & = & c_0^n \cdots  c_k^n c_{k+1}^{n-1}
 \cdots c_{n}^k  \\ & = & u^k ( c_0^{n-k} \cdots  c_k^{n-k} c_{k+1}^{n-(k+1)}
 \cdots c_{n-1}^1 ).
\end{eqnarray*}

For $i=1,\ldots, r$, let $k_i$ be the largest power of $u$ to divide all the 
elements
of the $i$th row of $\theta$. Set $P$, $Q$, $R$ to be diagonal matrices
whose $i$th diagonal entries are respectively
\begin{equation}
\begin{array}{c|c|c}
P_{ii} & Q_{ii} & R_{ii} \\ \hline
 c_0^{n-k} \cdots  c_k^{n-k} c_{k+1}^{n-(k+1)}
 \cdots c_{n-1}^1  & u^k & c_1^1c_2^2\cdots c_{k-1}^{k-1} c_k^k \cdots c_n^k \\
 \| & \| & \| \\
 e_1^-\cdots e_l^- & u^k & e_1^+\ldots e_k^+
\end{array}
\end{equation}
for $k=k_i$ and $l=n-k$. Then, 
\begin{lem} \label{lem:PQR}
\begin{eqnarray} P\cdot Q  & = & (\prod_{1}^{n}e_j^-) R.
\end{eqnarray}
\end{lem}

\section{Proper Transform of Sheaves} \label{sec:ProperTransform} 
Fix a point $x\in X$. Denote the blow-up of $X$ at $x$ by $\blow: \wtd{X} \lra X$
and the exceptional divisor by $\sE$.

\label{heckesect}
Let $T$ denote the zero dimensional subscheme corresponding to $x$. 

For a given coherent sheaf $\F$ on $X$, we call the pullback $\blow^* \F$ the {\em
total} transform of $\F$. Let 
$$\F^\sE:=\torsh{1}{\O_{\wtd{X}}}{\blow^*\F}{\O_{\wtd{X}}(\sE)_\sE}.$$ 
$\F^\sE$ coincides with the subsheaf of sections of $\blow^*\F$ supported along
$\sE$.

\begin{defn}
The {\em proper transform} of $\F$ is defined as the quotient  
$\blow^*\F/\F^\sE$ and will be denoted by $\wtd{\blow}\F$ , or simply $\wtd{\F}$ when
suitable.
\end{defn}
The following sequences are exact:
\begin{displaymath}
\begin{array}{lclclclcl}
0 & \lra & \O_{\wtd X} & \lra & \O_{\wtd X}(\sE) & \lra & \O_{\wtd X}(\sE)_\sE & \lra &0\\
0 & \lra & \F^\sE & \lra & \blow^*\F & \lra & \wtd \F & \lra &0 \\
0 & \lra & \wtd \F & \lra & \blow^*\F(\sE) & \lra & \blow^*\F(\sE)_\sE & \lra &0.
\end{array} 
\end{displaymath}
To see the exactness of the latter two, 
tensor  first sequence with $\blow^*\F$ and split the resulting sequence into two short exact sequences.

The definitions of proper transform for divisors and sheaves are compatible with each other:
For a given effective divisor $\sD$ on $X$, denote its ideal sheaf by $\J_\sD$, then 
$\wtd{\J}_\sD=\J_{\wtd\sD}$.

\begin{prop} \label{prop:MainProper}
Given 
a coherent sheaf $\F$ on $X$, 
\begin{enumerate}
\item If $\F$ is torsion-free, then $\F^\sE$ and $\wtd \F$ coincide with the torsion and torsion-free
parts of $\blow^*\F$ respectively.
\item Let $\J$ be the ideal sheaf of $x\in X$. Then, 
$\wtd{\J} = \O_{\wtd{X}}(-\sE)$.
\item If $\F$ is locally free at $x$, then $\blow^*\F=\wtd{\F}$.
In particular, if $\F$ is locally free, the conclusion holds.
\item If $x \notin \supp \F$, then $\blow^*\F=\wtd{\F}$ and 
$\blow_*\wtd \F=\blow_*\blow^*\F=\F$.
\item Given a locally free sheaf $L$ on $S$, then 
$(\F\otimes L)^\sE\isom \F^\sE\otimes L$ 
and $\wtd{(\F\otimes L)}\isom \wtd{\F} \otimes L$.
\end{enumerate}
\end{prop}
\begin{proof}
(1) follows as 
$\blow: \wtd{X}\backslash E \lra X \backslash \{ x \}$ is an isomorphism, 
$\blow^*\F$ is torsion-free over the open set $\wtd{X}\backslash E$ 
and the torsion locus is $\sE$.
(2) and (3) follow from (1).
\manudraftorfinal{}{and the fact if $\F_x$ is a free $\O_{X,x}$-module, 
then $\blow^*\F|E$ is locally free.} (4) is clear.
(5) follows from the locally freeness of the sheaf $L$. 
\end{proof}

\begin{lem} \label{lem:PushFwd}
Given a coherent sheaf $\M$ on $X$,
$$R^0\blow_*\blow^*\M=\M \; \mathrm{and} \; R^i\blow_*\blow^*\M=0 \; \textnormal{for all} \; i>0.$$
\end{lem}
\begin{proof}
The result holds for $\M=\O_X$. It holds for locally free sheaves by the projection
formula and for arbitrary coherent sheaves by the existence of 
locally free resolutions on smooth projective scheme $X$.
\end{proof}

\begin{lem} \label{lem:Purity}
Given a coherent sheaf $\M$ on $X$ with $\dimh(M_x)=1$,
\begin{enumerate}
\item
$\M^\sE\isom \O_\sE(-1)^{\oplus m}$, where $m=\dim_{k(x)} M_x\otimes k(x)$,
\item
$R^0\blow_*\wtd {\M}=\M$ and $ R^i\blow_*\wtd {\M}=0$ for all  $i>0.$
\item 
If $\M_x$ is torsion, then for all $y\in \sE$, $\dimh(\wtd{\M}_y)=1$,
\item 
$\torsh{i}{\O_{\wtd X}}{\wtd M}{\O_\sE}=0$ for all $i>0$.
\item $\extsh{0}{\O_{\wtd X}}{\O_\sE}{\wtd M}=0$.
\item If $\M$ is pure of dimension $1$, then $\sE\nsubseteq \supp \wtd{\M}$.
\end{enumerate}
\end{lem}

\newcommand{\N}{N}
\begin{assumption}
From now on, assume $\M$ is a coherent sheaf on $X$ with $\dimh(\M_x)=1$. 
Let $\N:=\ker(\M\lra \M_x)$.
\end{assumption}

\begin{lem} \label{lem:3rdCol}
Given an exact sequence of sheaves on $X$
$$0 \lra \N \lra \M \stackrel{\ev_T}{\lra} \M_T \lra 0$$
where $\ev_T$ is the evaluation map, then the sequence
$$0 \lra \wtd{N} \lra \blow^*\M \lra \blow^*\M_T \lra 0 $$
is exact and $\wtd N = \wtd\M(-\sE)$.

For any divisor $\sD$ with $mult_x \sD=1$,  $\wtd{M(-\sD)}=\wtd N(-\wtd \sD)$.

\end{lem}
\begin{proof}
The map $N\lra M$ is  $0$ at $x$ and an isomorphism of the fibers away from $x$.
Hence the map $\blow^*\N \lra \blow^*\M$  vanishes along $\sE$ and it is an isomorphism 
away from $\sE$. As a result, the kernel and the image 
of this map are $N^\sE$ and $\wtd N$. This proves 
the exactness of the above sequence.

The proper transform $\wtd\M$ fits into the exact sequence: 
$$ 0 \lra \wtd\M \lra \blow^*\M(\sE) \lra \blow^*\M(\sE)_\sE \lra 0.$$
Tensoring this sequence by $\O(-\sE)$ shows that $\wtd N = \wtd\M(-\sE)$.

Given such a divisor $\sD$, we see that $\blow^* =  \wtd \sD + \sE$. 
Starting from $\wtd N = \wtd\M(-\sE)$ and tensoring both sides by $\O(-\wtd \sD)$, 
we get 
$\wtd N(-\wtd \sD) = \wtd\M(-\blow^*)=\wtd{M(-\sD)}$.
\end{proof}

\begin{lem} \label{lem:ExactOnE}
The sequence $$0 \lra M^\sE \lra \blow^*\M_\sE \lra \wtd\M_\sE \lra 0 $$ is exact.
\end{lem}
\begin{proof} Let $K:=\ker(\blow^*\M_\sE \lra \wtd\M_\sE)$.
The following diagram is exact:

$$\xymatrix{ & & 0 \ar[d] & 0 \ar[d] & \\
 & & \M^\sE \ar@{=}[r] \ar[d] & K \ar[d] & \\
0 \ar [r]  & \wtd N \ar[r] \ar@{=}[d]& 
\blow^*\M \ar[r]  \ar[d] & \blow^*\M_\sE  \ar[r] \ar[d] & 0 \\
0 \ar [r]  & \wtd N \ar[r] & 
\wtd \M \ar[r]  \ar[d] & \wtd \M_\sE  \ar[r] \ar[d] & 0 \\
& & 0  & 0 & 
}.$$

\end{proof}

\subsection{Two Sheaves} \label{sec:TwoSheaves}
\renewcommand{\e}{\F}
If $\dimh(\M_x)=1$, then $\M$ has a two-step locally free resolution over an open subscheme
$U$ containing $x$:
$$ 0 \lra \F_1 \lra \F_0 \lra \M_U \lra 0.$$

Denote $\blow^{-1}(U)$ by $V$. This data fits into an exact diagram:

\begin{equation} \label{diag:onS} \tag{$\dagger$}
\xymatrix{ 
& & & 0 \ar[d] & \\
&  0 \ar[d] & & N_U \ar[d] & \\
0 \ar [r]  & \F_1 \ar[r]^\Phi 
\ar[d] & \F_0 \ar[r]  \ar@{=}[d] & M_U \ar[r] \ar[d] & 0 \\
0 \ar [r]  & {\mathcal H} 
\ar[r]  
\ar[d] & \F_0 \ar[r]  & 
M_T  \ar[r] \ar[d] & 0 \\
 & N_U \ar[d] & & 0 & \\
 &   0 &   &  & 
} 
\end{equation}
Here, $\sh$ and $N$ are the kernels of $\F_0\lra M_T$ and $M \lra M_T$ respectively.
The stalk $\M_x$ is torsion $\O_{X,x}$-module if and only if $\rank\F_1=\rank\F_0$.

\newcommand{\s}{{\mathcal S}}

\begin{lem} \label{lem:2ndRow}
Given an exact sequence of coherent sheaves on $U$
$$ 0 \lra \s \lra \F \lra \Q \lra 0 $$
with $\F$ locally free and $\Q$ torsion, 
then 
$$ 0 \lra \wtd\s \lra \blow^*\e \lra \blow^*\Q \lra 0 $$
is exact on $\wtd U$.
\end{lem}
\begin{proof}
The sequence 
$
\blow^*\s \lra \blow^*\F \lra \blow^*\Q \lra  0  
$
is exact. The first homomorphism factors through $\wtd\s$ because
$\wtd\s$ is the torsion-free quotient of $\blow^*\s$.  
The homomorphism $\wtd\s\lra \blow^*\F$ is generically injective since $\blow^*\Q$ is
torsion and injective since $\wtd \s$ is torsion free.
Consequently, the image of $\blow^*\s$ in $\blow^*\e$ coincides with $\wtd\s$. 
\end{proof}

\begin{defn}
Given a projective scheme $X$, a normal crossing divisor $\Sigma$, a locally free $\O_X$-module $\F$
and a locally-free $\O_\Sigma$-module $M$ together with a surjection $\phi: \F \lra M$, the coherent
sheaf $\ker\phi$
is called the Hecke transform of $\F$ with respect to $M$ and $\phi$.
\end{defn}

\begin{rk}
Hecke transforms are locally free sheaves on $X$. 
\end{rk}

\begin{prop}
\begin{enumerate}
\item If $\M_x$ is torsion $\O_{X,x}$-module, then $\sh$ 
is a torsion-free $\O_X$-module of the same
rank as $\e_1$.
\item The following diagram is exact:
\begin{equation} \label{diag:ProperS'} \tag{$\ddagger$}
\xymatrix{ 
& & & 0 \ar[d] & \\
&  0 \ar[d] & & \wtd  N_V \ar[d] & \\
0 \ar [r]  & \blow^*\F_1 \ar[r]^\Phi 
\ar[d] & \blow^*\F_0 \ar[r]  \ar@{=}[d] & \blow^*\M_V \ar[r] \ar[d] & 0 \\
0 \ar [r]  &\wtd  {\sh} 
\ar[r]  
\ar[d] & \blow^*\F_0 \ar[r]  & 
\blow^*\M_T  \ar[r] \ar[d] & 0 \\
 & \wtd  N_V \ar[d] & & 0 & \\
 &   0 &   &  & 
} 
\end{equation}
\item The proper transform $\wtd \sh$ of $\sh$ 
is a Hecke transform of $\e_1$
along $\sE$. In particular, $\wtd \sh$ is  locally free.
\item Given $\blow: {\wtd X} \lra X$ as before.
Then $\blow_* (\ref{diag:ProperS'}) = (\ref{diag:onS})$.

\end{enumerate}
\end{prop}

\begin{proof}
(1) The sheaf $\sh$ 
is torsion-free since any nontrivial subsheaf of a torsion-free
sheaf is torsion-free. The sheaves $\sh$ 
and $\F_0$ are of the same ranks since they are isomorphic away from $T$.

(2) The exactness of the second row and the third column follows from Lemmas 
\ref{lem:2ndRow} and \ref{lem:3rdCol}. The exactness of the first column is 
as a consequence. 

(3) The locus $\blow^{-1}(T)$ is a normal crossing divisor and $\blow^*\M_T$ a vector bundle 
on this divisor. The proper transform $\wtd  {\sh}$ of $\sh$
is the kernel of $\blow^*\F_0\lra \blow^*\M_T$ which proves it is a Hecke transform.

(4) follows from Lemmas \ref{lem:PushFwd} and \ref{lem:Purity}. 
\end{proof}

\begin{proof}[Proof of Lemma \ref{lem:Purity}]
(1) The stalk $\M_x$ has a two-step resolution by free $\O_{X,x}$-modules.
Therefore, there exists an open neighborhood $U$ of $x$ on which $\M_U$
has a locally free resolution
 $$ 0 \lra \F_1 \lra \F_0 \lra \M_U \lra 0.$$
The locally free sheaves $\F_i$ are of the same rank, say $r$. Let $m$ be the fiber
dimension of $M$ at $x$, then $m\leq r$. Denote $\blow^{-1}(U)$ by $V$. Because
$\sE\subset V$, $\M^\sE=\torsh{1}{\O_{V}}{\blow^*\M_U}{\O_{V}(\sE)_\sE}$. Using the locally free
resolution of $M_U$, we see that $\M^\sE\isom \O_\sE(-1)^{\oplus m}$ where
$m=\dim_{k(x)} M_x \otimes k(x)$.  

(2) 
The sequence $
\begin{array}{lclclclcl}
0 & \lra & \F^\sE & \lra & \blow^*\F & \lra & \wtd \F & \lra &0 
\end{array}$
is exact. The first term is isomorphic to $\O_\sE(-1)^{\oplus m}$. Therefore
$R^i\blow_*\M^\sE=0$ for all $i$, $R^0\blow_*\wtd {\M}=\M$ and 
$ R^i\blow_*\wtd {\M}=0$ 
for all  $i>0.$

(3)
For all $y\in \sE$, $\dimh(\blow^*\M_y)=1$ and $M_y$ is a torsion $O_{\wtd X, y}$-module.
The sheaf $\wtd N$ is a subsheaf of $\blow^*\M$ and $\wtd\M = \wtd N(\sE)$. By Cor.
\ref{cor:AusBuchs},
$\wtd\M_y$ has the same properties.

(4)
The sheaf $\wtd\M_V$ has a two-step locally free resolution
$$ 0 \lra \F_1(\sE) \lra \wtd \sh(\sE)\lra \wtd\M_V \lra 0,$$
where $\wtd \sh$ is defined by a Hecke transform as $\ker(\F_0 \lra \blow^*\M_\sE)$ and
hence locally free. Consequently, $\torsh{i}{\O_{\wtd X}}{\wtd
M}{\O_\sE}= 0$ for $i\geq 2$.
 
Apply $\O_\sE\otimes_{\O_{\wtd X}} \bullet $ to the exact sequence
$$0 \lra \M^\sE \lra \blow^*\M \lra \wtd\M \lra 0.$$
The sheaf $\M^\sE\isom \O_\sE(-1)^{\oplus m}$. Hence 
$\torsh{1}{\O_{\wtd X}}{\M^\sE}{\O_\sE}\isom \O_\sE^{\oplus m}$. Similarly, 
$\torsh{1}{\O_{\wtd X}}{\blow^*\M}{\O_\sE}=\M^\sE(-\sE)\isom \O_\sE^{\oplus m}$.
By Lemma \ref{lem:ExactOnE}, 
$\torsh{1}{\O_{\wtd X}}{\wtd\M}{\O_\sE}=0$.

(5) Apply $\extsh{}{\O_{\wtd X}}{\bullet}{\wtd\M}$ to 
the exact sequence 
$$ 0 \lra \O_{\wtd X}(-\sE) \lra \O_{\wtd X} \lra \O_\sE \lra 0.$$
The result is 
$$ 0 \lra \extsh{0}{\O_{\wtd X}}{\O_\sE}{\wtd\M} 
\lra \wtd\M \lra \wtd\M(\sE) \lra  \wtd\M(\sE)_\sE \lra 0.$$
Then, 
$\extsh{0}{\O_{\wtd X}}{\O_\sE}{\wtd\M}=\torsh{1}{\O_{\wtd
X}}{\wtd\M}{\O_\sE}=0$. The latter sheaf is trivial by (4).

(6) If $\M$ is pure of dimension $1$, then $\wtd \M$ is torsion and $\dimh(\wtd \M_y)=1$
for all $y\in \wtd X$. Consequently, $\wtd\M$ is pure of dimension $1$. As $\sE$ is not
an associated point of $\wtd \M$ by (4), $\sE\nsubseteq \supp \wtd \M$.
\end{proof}

\subsubsection{Morphisms} Let $\M_1, \M_2$ be coherent sheave on $X$. Then a
homomorphism $\phi: \M_1 \lra \M_2$ induces $\phi^*:\blow^*\M_1\lra\blow^*\M_2$ 
with $\blow^*(\M_1^\sE)\subset \M_2^\sE$ and denote the induced morphism on the
quotients
$\wtd \M_1 \lra \wtd \M_2$ by $\wtd \phi$. 
For coherence, we also use $ \wtd{\blow}\phi$ for  $\wtd\phi$.

\begin{assumption}
$\M_1, \M_2$ are pure of dimension $1$.
\end{assumption}

\begin{lem}\label{lem:qi}
For any homomorphism $\phi: \M_1\lra \M_2$, $\blow_*\wtd {\blow}\phi=\phi$. 
For any homomorphism $\psi: \wtd\M_1\lra \wtd\M_2$, 
$\wtd{(\blow_*\psi)}=\psi$. 
\end{lem}
\begin{proof}
The kernel of $\phi-\blow_*\wtd \phi$ is $0$-dimensional. Because $\M_1$ is pure 
of dimension $1$, the kernel is trivial and hence $\phi=\blow_*\wtd \phi$.
The kernel of $\psi-\wtd{(\blow_*\psi)}$ is contained in the zero dimensional 
subscheme $\sE\cap \supp \wtd \M_1$. Hence, $\psi=\wtd{(\blow_*\psi)}$.
\end{proof}
\begin{lem} \label{lem:Injectivity}
A homomorphism $\phi: \M_1 \lra \M_2$ is injective if and only if 
$\wtd \phi$ is injective. In this case, $\blow_*(\coker\wtd\phi)=\coker\phi$.
\end{lem}
\begin{proof}
($\Leftarrow$) As $M_i$ are pure of dimension $1$, $M_i=\wtd \M_i$ and 
$\phi=\blow_*\wtd {\blow}\phi$. The injectivity follows from the left exactness of $\blow_*$.

\noindent
($\Rightarrow$) $\wtd \phi$ is injective away from $\sE$, thus $\ker\wtd\phi\subset
E\cap \supp \wtd \M_1$, hence trivial.   
\end{proof}

\begin{lem}\label{lem:EPpres}
Let $\M$ be a pure sheaf of dimension $1$. Then,
\begin{eqnarray*}
\chi(\blow^*\M) & = & \chi(\M), \\
\chi(\wtd \blow \M) & = & \chi(\M). 
\end{eqnarray*}
\end{lem}
\begin{proof}
Lemma \ref{lem:PushFwd} (resp. Lemma \ref{lem:Purity}) shows 
that the sheaf cohomology of $\blow^*\M$ (resp. $\wtd \blow \M$) match the sheaf
cohomology of $M$, hence $\chi(\blow^*\M)  = \chi(\M) = \chi(\wtd \blow \M)$.
\end{proof}

\subsection{Parabolic Case} \label{subsec:ParPropTr}

Fix a divisor $\sD$ on $X$.

\begin{defn}
Given a parabolic sheaf
$\M_\bullet$ of dimension $1$ on $X$ with parabolic divisor $\sD$, 
the {\em proper transform} of $\M_\bullet$ is defined by setting 
$\wtd \blow (\M)_\alpha:=\wtd\blow(\M_\alpha)$ for $\alpha \in \RR$.
Let $\wtd \blow M_\bullet$ have the same  weights as $\M_\bullet$.
\end{defn}

\begin{defn}
Given a parabolic sheaf
$\wtd \M_\bullet$ of dimension $1$ with parabolic divisor $\blow^*\sD$,
define the {\em push-forward} of $\wtd \M_\bullet$ by setting
$\blow_* (\wtd \M)_\alpha:=\blow_*(\wtd \M_\alpha)$ for  $\alpha \in \RR$.
Let $\blow_*\wtd  M_\bullet$ have the same  weights as $\wtd \M_\bullet$.
\end{defn}

\begin{prop}
The proper transform $\wtd \blow \M_\bullet$
of  $\M_\bullet$ with divisor $\sD$  is a parabolic sheaf
with divisor $\blow^*\sD$, $\blow_*\gr{F}{i}\wtd \blow \M=\gr{F}{i}M$ 
for all $i$ and $\chip(\wtd \blow \M_\bullet)  = \chip(\M_\bullet)$.
\end{prop}

\begin{proof}
We give a proof for the proper transform. First, $\wtd \blow (\M)_0 = \wtd\blow \M$.
By assumption, $\M_\bullet$ is parabolic. The sheaf
$\M_\beta$ is a subsheaf of $\M_\alpha$  for $\beta \geq \alpha$.
The previous lemma implies $\wtd \M_\beta$ is a subsheaf of $\wtd \M_\alpha$
and $\blow_*\gr{F}{i}\wtd \blow \M=\gr{F}{i}M$ 
for all $i$. 
In addition, 
$$\wtd \blow (\M)_{\alpha+1}=\wtd \blow (\M_\alpha(-\sD))=
\wtd \blow (\M_\alpha)(-\sD))=\wtd\M_\alpha(-\blow^*\sD).$$
The weights of $\wtd \blow \M_\bullet$ are the weights of $\M_\bullet$. 
These make $\wtd \blow \M_\bullet$ into a 1-dimensional parabolic sheaf.

The parabolic Euler characteristic is preserved 
since 
$$\wtd{M(-\sD)}=\wtd \M(-\blow^*\sD)$$ and $\blow_*\gr{F}{i}\wtd \M=\gr{F}{i}M$ 
for all $i$:
$$
\begin{array}{lclcl}
\chip(\wtd\blow \M_\bullet) & = & \chi(\wtd\blow \M(-\blow^*\sD)) & +  
& \sum \alpha_i \chi(\gr{F}{i}\wtd\blow \M)\\
& = & \chi(\M(-\sD)) & + &  \sum \alpha_i \chi(\gr{F}{i} \M)\\
& = & \chip(\M_\bullet).& &
\end{array}
$$
\end{proof}

\section{Addition and Deletion} \label{sec:AddDelete}
In this section, we introduce the addition and deletion operations
for parabolic sheaves.

\subsubsection{Deletion}
Let $\P_\bullet$ be a parabolic sheaf on $X$ with divisor $\sD$ 
and $\sD',\sE$ be a effective Cartier divisors in $X$ 
such that $\sD=\sD'+\sE$. Because $\dim(\supp(\P)\cap\sD)<\dim(\supp(\P))$, 
the same holds for $\sE$ and $\sD'$ as well. We set $\P':=\P(-\sE)$.
One can put a parabolic structure on the sheaf $\P'$ whose parabolic divisor is $\sD'$: 
For $0 \leq \alpha < 1$,
set $\P'_\alpha:=\P'\cap \P_\alpha$. Extend this to a parabolic structure by 
setting $\P'_\alpha:=\P'_{\alpha-\floor{\alpha}}(\floor{\alpha}\sD')$. 
We call $\P'_\bullet$ the {\em deletion of $\sE$ from the divisor of  $\P_\bullet$}, and denote it 
by $\del{\sE}(\P_\bullet)$.

\subsubsection{Addition}\label{sssec:add}
Conversely, given a parabolic sheaf $\P'_\bullet$ with divisor $\sD'$ and an effective
divisor $\sE$ such that $\dim(\supp(\P')\cap\sE)<\dim(\supp(\P'))$, 
one can put a parabolic structure with parabolic divisor $\sD=\sD'+\sE$ on
$\P=\P'(\sE)$ by setting $\P_0 = \P$ and $\P_\alpha=\P'_\alpha$ for $0<\alpha< 1$
and extending this to $\RR$ in the usual way. 
Denote $\P_\bullet$ by $\add{\sE}(\P'_\bullet)$, and call it the {\em addition of $\sE$ to the divisor of  $\P'_\bullet$}. 


\begin{assn}\label{assn:second}
One has either
$$\alpha_0(\P)=0 \; \textnormal{ and } \; F_1\P=\P(-\sE),$$
or 
$$\alpha_0(\P)>0 \; \textnormal{ and } \; \supp(\P)\cap\sE=\emptyset .$$
\end{assn}

\begin{lem} \label{lem:preserveSecAssn}
Let $\P_\bullet$ be a parabolic sheaf on $X$ satisfying Condition \ref{assn:second}, and 
$\S_\bullet$ be a saturated subsheaf endowed with the induced parabolic structure. Then $\S_\bullet$ 
also satisfies Condition \ref{assn:second}.
\end{lem}

\begin{proof}
The induced parabolic structure is defined by the formula $\S_\alpha =\S\cap\P_\alpha$ for 
all $0\leq \alpha <1$. For small enough $\varepsilon >0$, one then has 
$\S_\varepsilon=\S\cap\P_\varepsilon=\S\cap\P(-\sE)$. Because $\S\subset\P$, this then implies 
$\S_\varepsilon=\S(-\sE)$. Assuming $\supp(\S)\cap\sE\neq\emptyset$, one has $\S_\varepsilon\neq\S$. 
Letting $\varepsilon\to 0$, this implies $\alpha_0(\S)=0$, and $F_1\S=\S(-\sE)$ as claimed. 
\end{proof}

\begin{prop}\label{prop:adddelinv}
For any parabolic sheaf $\P'_\bullet$ with divisor $\sD'$ and any $\sE$ satisfying the assumption of Subsection 
\ref{sssec:add}, one has $\del{\sE}(\add{\sE}(\P'_\bullet))=\P'_\bullet$. 
For any parabolic sheaf $\P_\bullet$ and $\sE$ satisfying Condition \ref{assn:second}, 
one has $\add{\sE}(\del{\sE}(\P_\bullet))=\P_\bullet$. Finally, 
$\chip(\del{\sE}(\P_\bullet))=\chip(\P_\bullet)$.
\end{prop}

\begin{proof}
The first two statements are clear from the definitions. For the third statement, notice that 
$\P(-\sD)=(\P(-\sE))(-\sD')$ and that for all $0<\alpha\leq 1$, 
one has $\P_\alpha=\del{\sE}(\P_\bullet)_\alpha$, in particular 
$\gr{F}{\alpha}(\P_\bullet)=\gr{F}{\alpha}(\del{\sE}(\P_\bullet))$. 
\end{proof}

\section{Spectral Triples}\label{sec:SpecTr}
\newcommand{\X}{x}
\newcommand{\Y}{y}

From now on, let $C$ be a smooth projective curve over the field $\kk$ and 
$\sP$ an effective Cartier divisor on $C$.

\begin{defn}
A spectral triple $(Z,\Sigma, M)$ consists of a smooth  surface $Z$ together with a morphism $\blow: Z \lra C$, an effective divisor $\Sigma \subset Z$ and
a rank one torsion free $\O_\Sigma$-module $M$  so that 
\begin{tabular}{lcl}
(1) & $\pi$ & $: Z \lra C$ is flat, \\
(2) & $\pi_{|\Sigma}$ & $: \Sigma \lra C$ is finite and flat.
\end{tabular}
\end{defn}

\begin{defn}
A parabolic spectral triple $(Z,\Sigma, M_\bullet)$ is a spectral triple so that
$M_\bullet$ is a parabolic sheaf whose parabolic divisor  is $\blow^*(\sP)$.
\end{defn}

\begin{rk}
In fact, it would be enough to consider the pair $(Z,M)$ as $\Sigma$ is determined
by $\Sigma=\supp(M)$. We include $\Sigma$ in the definition for expositional purposes.
\end{rk}

Starting from such a Higgs bundle $(\E,\theta)$ one can define a spectral 
triple $(Z^\oo,\Sigma^\oo, M^\oo)$:
Set $C^\oo:=C -\sP$ and $\AA^1:=\spec(k[\lambda])$. Here,

\begin{align}
Z^\oo & :=  C^\oo\times\AA^1 \;\text{ with the obvious choice for}\;
\blow_\oo: Z^\oo \lra C, \\
\Sigma^\oo &:=\div{\det(\lambda\, \Id_\E - \theta)} \;\text{and}\\
M^\oo &:= \coker(\lambda \,\Id_\E - \theta).
\end{align}

\begin{defn}
A compactification of the triple $(Z^\oo,\Sigma^\oo, M^\oo)$ is a spectral
triple $(Z,\Sigma, M)$ and an open immersion $i: Z^\oo \lra Z$ so that
\begin{enumerate}
\item $Z$ is a connected smooth projective surface, 
\item $i(\Sigma^\oo)$ is contained in $\Sigma$ as a dense open subset,
\item $i^*\M=M^\oo$, 
\end{enumerate}
so that the following diagram is commutative:
\[
\xymatrix{
Z^\oo \ar[dr]_{\pi_\oo} \ar[rr]^{i} &  & Z \ar[dl]^{\pi}\\
& C. &
}
\]
\end{defn}

By (1) and (2), $Z^\oo$ is dense inside $Z$ and
 $\Sigma$ is the closure of $i(\Sigma^\oo)$ in $Z$.


There are many compactifications of the triple
$(Z^\oo,\Sigma^\oo, M^\oo)$: Given one compactification, one can obtain other
via suitable blow-ups.

The standard choice for compactification of $(Z^\oo,\Sigma^\oo, M^\oo)$
is 

\begin{tabular}{llcl}
 & $Z^\sP$ & $:=$ & $\PP_C(\O\oplus \O(-\sP))$,\\
 & $\Sigma^\sP$ & $:=$ & $\div{\det(\X_\sP \, \Id_\E - \Y_\sP \, \theta)}$ and \\
 & $M^\sP$ & $:=$ & $\coker(\X_\sP \,\Id_\E - \Y_\sP \, \theta)$.
\end{tabular}

\begin{rk}
 $\Sigma^\sP$ does not meet the line at infinity $(\Y_\sP)$.
\end{rk}

\subsection{Naive Compactification: The Surface and The Curve} 
We describe a compactification with $Z^0=C\times \PP^1$. The spectral 
curve $\Sigma^0$ is determined by $\Sigma^\oo$ as its closure. The sheaf $M^0$ is 
the cokernel of a morphism between locally free sheaves as was $M^\sP$.

\begin{rk}
$\Sigma^0$ meets the line at infinity $(\Y_0)$
over the points $p\in\sP$ for which $\theta_p$ is not a nilpotent endomorphism. 
\end{rk}


\subsection{Naive compactification: The Sheaf}
Fix a section $s$ so that $\sP=\div{s}$ and let  
$\Theta_0 := s \, \X_0\, \Id_\E - \Y_0 \, \theta: \E \lra \E(\sP)\otimes\O_{Z^0}(1)$.
Then $$\div{\det\Theta_0}=\Sigma^0+\div{\det\coker(\theta)}.$$
We find another map $\Phi_0$ related to $\Theta_0: \E \lra \E(\sP)\otimes\O_{Z^0}(1)$ so that 
$$\Sigma^0=\div{\det \Phi_0}.$$ 
Recall that $\F:=\ker(\E(\sP) \lra \coker\, \theta(-\sP)_\sP)$ and $\F
\hookrightarrow \E$ is denoted by $Q$. 
The map $\theta: \E\lra \E(\sP)$ naturally factors through $\F(\sP)$.
Similarly, $\Theta_0$ factors through $\F(\sP)\otimes \O_{Z^0}(1)$. Denote the resulting
map by $\Phi_0$. The following diagram is exact. From this diagram, we see that
$$\supp M^0=\div{\det \Phi_0}=\Sigma^0.$$

\begin{eqnarray} \label{diag:Naive3}
\xymatrix{ 
& & & 0 \ar[d] & \\
&  0 \ar[d] & & M^0 \ar[d] & \\
0 \ar [r]  & \E \ar[r]^{\Theta_0} \ar[d]_{\Phi_0} & \E(\sP)\otimes\O_{Z^0}(1)\ar[r]  \ar@{=}[d] & \pmapz^*\coker\,\theta \ar[r] \ar[d] & 0 \\
0 \ar [r]  & \F(\sP)\otimes\O_{Z^0}(1) \ar[r]^{Q \otimes 1}  \ar[d] & \E(\sP)\otimes\O_{Z^0}(1)  \ar[r]  & \pmapz^*\coker\,\theta_{D} \ar[r] \ar[d] & 0 \\
 & M^0 \ar[d] & & 0 & \\
  &   0 &   &  & 
} & & 
\end{eqnarray}

\subsection{}
We use the ideas developed in
Section \ref{sec:ProperTransform} to compare $M^0$ and $M^\sP$ using the diagram 
$$
\xymatrix{
& \ar[ld]_{\pmapz} Z  \ar[rd] ^{\pmapP} & \\
Z^0  \ar@{-->}[rr] && Z^\sP.}
$$

Set
\begin{itemize}
\item $T^-:=\Sigma^-\cap (y_0)$
and 
$T^+:=\Sigma^+\cap (x_\sP)$,
\item 
$\F_1^0:=\E$, $\F_0^0:=\F(\sP)\otimes\O_{Z^0}(1)$
and 
$\F_1^\sP:=\E$, $\F_0^\sP:=\E(\sP)\otimes\O_{Z^\sP}(1)$.
\end{itemize}
View $\F_1^\sB\lra \F_0^\sB$ as a two-step locally
free resolution for the sheaf $M^\sB$ on $U=Z^\sB$. Set 
$$\sh_\sB:=\ker(\F_0^\sB \lra \M^\sB_{T_\sB}).$$

\begin{prop} \label{prop:Relation}
For $\sB=0,\sP$,  $T^\sB$ is  linear, \\$Z=\fbl{T^\sB}Z^\sB$ and 
 $\blow_\sB=\blow_{T^\sB}$,
\begin{eqnarray}
\wpmapP\sh^\sP & = & \wpmapz\sh^0 \\
\wpmapP N^\sP & = & \wpmapz N^0 \\
 \wtd {\Sigma^0} & = &\wtd{\Sigma^\sP}.\label{prtrspectralcurve}
\end{eqnarray}

\end{prop}

We denote the proper transformed spectral curve (\ref{prtrspectralcurve}) by $\Sigma$.

\begin{proof}
For simplicity, assume $\sP=n\cdot \textnormal{pt}$ with local coordinate $u$.
at $\textnormal{pt}$.
Subscheme $T^\sP$ is linear because $T^\sP=(u^n)\cap (t)=(u^n,t)$ where $t$ is $x_\sP$ or
$y_D$. The pair $u$ and $t$ is clearly transversal. Rest of (1) is clear.
Lemma \ref{lem:PQR} says that the following is diagram commutative and 
$\wtd \pmapP\sh^\sP  = \wtd
\pmapz\sh^0$:
\begin{eqnarray*}
\xymatrix{
\wtd{\sh_\sP} \ar@{-->}[rr]^{(\isom)} \ar[d]_{P} & &  
\wtd{\sh_0} \ar[d]^R \\
\F(\sP)\otimes\O_{Z^0}(1) \ar[dr]_Q & & \E(\sP)\otimes\O_{Z^\sP}(1) \ar[dl]^{(\prod_1^n e_i^-)} \\
& \E(\sP)\otimes \O_{Z^\sP}(1) & 
}
\end{eqnarray*}
The rest follows.
\end{proof}

\newcommand{\absolute}{intermediate }

\section{Algebraic Nahm Transform} \label{trsect}
Let $(Z^\sP,\Sigma^\sP,M^\sP_{\bullet})$ be the standard parabolic spectral triple of the parabolic Higgs bundle
$(\E_{\bullet}, \theta )$ whose parabolic divisor $\sD$ is $\sP+\infty$. Define the following $0$-dimensional subschemes of $Z^\sP$: 
$$
    T^{+}:=\pi^*(\sP)\cap (x_\sP),
$$ 
the intersection of the pullback divisor $\pi^*(\sP)$ and the $0$-section 
in $Z^\sP$ and 
$$
    \Tt^{-}:=\pi^*(\infty)\cap  \Sigma^\sP, 
$$ 
the fiber of $\Sigma^\sP$ over $\infty \in Z^\sP$. 
Applying the ideas of  Section \ref{heckesect}  to the zero-dimensional
subscheme $T^{+}$ produces
a new spectral triple $(Z^\sP,\Sigma^\sP,N^\sP)$ out of $\M^\sP_\bullet$: 
$$ N^\sP:=\ker(M^\sP \ra M^\sP_{T^+}). $$
By definition, $N^\sP$ consists of the local sections of $M^\sP$ 
vanishing in $T^+$.

Next, define another spectral triple $(Z,\Sigma,N_{\bullet})$: 
The surface $Z$ is the blow-up $\pmapP$ of $Z^\sP$ at  $T^{+}$, the divisor 
$\sE^{+}$ is the exceptional divisor of $\pmapP$ and the coherent sheaf $N$ 
is defined as the proper transform  $\pmap_\wtd{\sP}(N^\sP)$ of $N^\sP$. 
The support $\Sigma$ of $N$ is the proper transform $\wpmapP(\Sigma^\sP)$ of 
$\Sigma^\sP$. Let $$N_{\bullet}:=\del{\sE^{+}}\wpmapP M^\sP_{\bullet}.$$  
The sheaf $N$ now has a parabolic structure whose divisor is 
$\wpmapP\sD=\pmapP^{\ast}\sD-\sE^{+}$. 
By our convention, the parabolic structure of $N^\sP_{\bullet}$ (as well as that of all the other sheaves involved further in the 
construction) has weights between $0$ and $1$ in all parabolic points. 

\begin{prop}\label{prop:condts}
If Condition \ref{assn:Main} holds for $(\E ,\theta )$, then Condition \ref{assn:second} holds for the 
parabolic sheaf $N_{\bullet}$ and the divisor $\sE^{+}$. 
\end{prop}
\begin{proof}
Indeed, Condition \ref{assn:second} says that in a
parabolic point $p\in\sP$ such that $\theta_p$ has a $0$ eigenvalue, the relation 
$\F=F_1 \E$ holds, and the smallest parabolic weight is $0$. 
This then implies that $\gr{F}{0}(\E)=\coker (\theta_p)$, and that the parabolic weight associated to this graded 
piece is $0$. By the definition of the parabolic structure of $M^\sP$, this is equivalent to saying that 
$\gr{F}{0}(M^\sP)=\coker (\theta_p)$ on the fiber $\sF=\pi^{-1}(p)$, with $0$ parabolic weight associated to this graded. 
Since the sheaf $\coker (\theta_p)$ is supported in the point $t=\sF\cap (x_\sP)$, we see that the 
$0$-weight subspace of $M^\sP|_{\sF}$ is precisely $M^\sP|_t$. Therefore, we have $F_1M^\sP=ker(M^\sP\to M^\sP|_t)$. 
Now let us blow up the point $t$, and call $\sE$ the exceptional divisor and $M$ the proper transform of $M^\sP$. 
Then for the parabolic structure of $M$ the relation $F_1M=ker(M\to M|_\sE)=M(-\sE)$ holds, and the corresponding 
parabolic weight is $0$. This shows that Condition \ref{assn:second} is true for $M$. It then follows for $N$ 
as well because of our convention of keeping the same weights for kernel and cokernel sheaves. 
\end{proof}

For simplicity,  write 
$\pi$ for the projection $\pi \circ \pmapP: Z^\sP\lra \CP$, whenever this does not create any
ambiguity, and do the same for all other projections to $\CP$ composed with blow-up maps. Similarly, identify the 
zero-dimensional
subcheme $T\subset Z^2$ with $\blow^{-1}(T)$ 
if $\blow:Z^1\lra Z^2$ does not affect $T$. 
For example, view $\Tt^{-}=\pi^*(\infty)\cap\Sigma^\sP$ both
as a zero-dimensional  subscheme  
of  $Z^\sP$ and  $Z$. 

Now, apply the blow-up construction of Section \ref{heckesect} 
to the zero-dimensional scheme $\Tt^{-}$ in $Z$, and obtain a parabolic 
spectral triple $(Z^{\abs},\Sigma^{\abs},N^{\abs}_{\bullet})$ called the 
\emph{\absolute parabolic spectral triple} of 
$(\E_{\bullet}, \theta )$: the surface $Z^{\abs}$ is the blow-up $\qmapP$ 
of $Z$ along $\Tt^{-}$, the coherent sheaf
$N^{\abs}$ is defined as the proper transform $\wqmapP\wpmapP N$ of $N$ 
with respect to $\pmapP\qmapP$. It is 
supported on the proper transform
$\Sigma^{\abs}=\wqmapP\Sigma$ of $\Sigma$ with respect to $\qmapP$.  Set 
$$N^{\abs}_{\bullet}:=\del{\sE^{+}}\wtd{(\pmapP\circ \qmapP)}(M^\sP_{\bullet}).$$
The parabolic divisor of $N^{\abs}_\bullet$ is $\sD^{\abs}=
\qmapP^*\wpmapP{\sD}$. Call $\widehat{\sE}^{-}$ the exceptional divisor of $\qmapP$. 
Set $\sP^{\abs}=\sD^{\abs}\setminus \pi^{-1}\infty$. 
We call $Z^{\abs}$ the \emph{\absolute surface}, 
$\Sigma^{\abs}$ the \emph{\absolute spectral curve} and $N^{\abs}$ the \emph{\absolute spectral sheaf}. 

It is possible to reconstruct the original parabolic Higgs bundle from the \absolute parabolic
spectral triple: by Lemma \ref{lem:Purity}, we have 
$$N^\sP=(\pmapP \circ \qmapP)_{\ast} N^{\abs} \; \mbox{and} \; 
M^\sP(-\sP)_\bullet=(\pmapP \circ \qmapP)_{\ast} N^{\abs}_\bullet(-\sP^{\abs}).$$
The parabolic vector bundle $\E_{\bullet}$ is
$\pi_{\ast}M^\sP(-\sP)_{\bullet}$, and the Higgs field $\theta$ 
is the direct image of multiplication map by the
global section $x_\sP$ on $Z^\sP$. 

The dual divisor $\widehat \sP$ is defined as the image of $\widehat T^-$ 
under $\widehat \pi: \PP^1\times \widehat \PP^1 \lra \widehat \PP^1$, i.e.
$$\widehat \sP:=\widehat \pi(\widehat T^-).$$
We equally set $\widehat \sD=\widehat \sP + \ti$.
All the maps we have constructed so far fit into the commutative diagram 
\begin{align}\label{absdiag}
\xymatrix{
       & & Z^{\abs} \ar[dl]_{\qmapP} \ar[dr]^{\qmapPt} & & \\
       & Z \ar[dl]_{\pmapP} \ar[dr]^{\pmapz} & & 
       \Zt \ar[dl]_{\pmaptz} \ar[dr]^{\pmaptP} & \\
       Z^\sP & & \CP \times \CPt & & \Zt^\Pt.}
\end{align}
Here, the surfaces $\widehat Z^{\widehat \sP}$, $\widehat Z$ and the related maps are defined in an analogous manner to above.
More precisely, recall that $\pmapP$ is the blow-up of $Z^\sP$ at the points $T^{+}$,
and $\pmapz$ is the blow-down of $Z$ along the proper transforms $\sE^{-}=\wpmapP(\pi^*(\sP))$ in $Z$ of the
fibers of $\pi$ in $Z^\sP$ over the points of $\sP$. As usual, call $\pi$ and $\pit$ the two projections of $\CP \times \CPt$. 
For any $p\in \sP$ the proper transform $\wpmapP(\pi^*(p))$ of the fiber over $p$ contracts into the 
point in $\CP \times \CPt$ which is the intersection of the infinity-fiber of $\pit$ with the fiber of $\pi$ over
$p$: denote by $T^{-}\subset\pit^*(\widehat{\infty})$ the union of these points for all $p\in \sP$; this is a finite set. 
Furthermore, recall that $\Tt^{-}$ is the intersection of $\Sigma^{0}$ and the infinity-fiber of $\pi$, or said differently, 
the intersection of the fibers of $\pit$ over $\Pt$ in $\CP \times \CPt$ and the infinity-fiber of $\pi$. 
Then, the map $\pmaptz$ is the blow-up of $\CP \times \CPt$ in the points $\Tt^{-}$, $\qmapPt$ is the blow-up of $\Zt$
in the points $T^{-}$, and finally $\pmaptP$ is the blow-down of the proper transform 
$\widehat{\sE}^{+}=\wpmaptz(\pit^*(\Pt))$ of the fibers of $\pit$ over $\Pt$ with respect to $\pmaptz$. 
Call $\Tt^{+}$ the finite set where these fibers contract in the $0$-section of $\pit$ in $\Zt^{\Pt}$. 
In other words, the relation between $\CP \times \CPt$, $\Zt^{\Pt}$ and $\Zt$ with respect to the points $\Pt$ 
and the projection $\pit$ is the same as the relation between $\CP \times \CPt$, $Z^\sP$ and $Z$ with respect to
the points $\sP$ and the projection $\pi$; whereas $Z^{\abs}$ is the fibered product of $Z$ and $\Zt$ over $\CP \times \CPt$. 
Therefore, $Z^{\abs}$ has two projections to projective lines: 
$\pi = \pi \circ \pmapP \circ \qmapP$ to the base $\CP$ of the
geometrically ruled surface $Z^\sP$, and $\pit = \pit \circ \pmaptP \circ \qmapPt$ to the base $\CPt$ of $\Zt^{\Pt}$. 
Let $\Mt^{\Pt}$ be the direct image sheaf $(\pmaptP \circ \qmapPt)_{\ast} (\add{\widehat{\sE}^{+}} N^{\abs})$ and denote by 
$\widehat{\Sigma}^{\Pt}$ its support. 

\begin{defn}
The direct image parabolic sheaf $\pit_{*} \Mt^\Pt_{\bullet}(-\Pt)$ on
$\CPt$ with parabolic points $\Dt$ will be called $\Et^{\prime}_{\bullet}$.
\end{defn}
By the definition of $\Mt^\Pt_{\bullet}$, $\Et^{\prime}_{\bullet}$ is
isomorphic to
$$
    \left(\pit_{*}\circ (\pmaptP)_{*} \circ (\qmapPt)_{*} \circ
\add{\hat{E}^{+}} \del{E^{+}}
     \circ \wqmapP \circ \wpmapP (M^{\sP}_{\bullet})\right)(-\Pt).
$$

Remark that by construction, the parabolic weights of $\Et^{\prime}_{\bullet}$ are between $0$ and $1$. 
Furthermore, similarly to $(x_\sP,y_\sP)$ on $Z^\sP$, there exists a pair of globally well-defined parameters
$(\Xt_{\Pt},\Yt_{\Pt})$ on $\Zt^\Pt$. 
\begin{defn}
Denote the direct image of multiplication by the global section $-\Xt_\Pt$ on $\Mt^\Pt_{\bullet}(-\Pt)$ by $\tt^{\prime}$. 
\end{defn}

\begin{rk}
One checks without difficulty that $\tt^{\prime}$ respects the parabolic filtration of $\Et^{\prime}_{\bullet}$, 
hence $(\Et^{\prime}_{\bullet},\tt^{\prime})$ is a parabolic Higgs bundle. Using the notation introduced in Section 
\ref{sec:outline}, the definitions above can be written 
$$
(\Et^{\prime}_{\bullet},\tt^{\prime})=\pi_H(\Mt^\Pt_{\bullet},-\Xt_\Pt), 
$$ 
or equivalently by Lemma \ref{lem:(-1)}
$$
(\Et^{\prime}_{\bullet},\tt^{\prime})=\pi_H((-1)^*_{\Zt^\Pt}\Mt^\Pt_{\bullet},\Xt_\Pt)
$$
\end{rk}

Notice that the construction of $(\Et^{\prime}_{\bullet},\tt^{\prime})$ only assumes that $\theta$ has first-order poles 
at finite distance and no singularity at infinity, but no assumption is made on the residues of $\theta$ in these 
singularities, nor about stability or the degree of $\E$. However, the reason why we introduced 
this construction is that under the assumptions of \cite{Sz}, the two definitions of Nahm transform agree: 
\begin{thm}\label{trthm}
Assume $(\E_{\bullet},\theta)$ satisfies the conditions of Section \ref{nahmsect}. Then, the parabolic Higgs bundles 
$(\Et^{\prime}_{\bullet},\tt^{\prime})$ and $(\Et^{tr}_{\bullet}, \tt )$ are isomorphic. 
\end{thm}
\begin{defn}
In the sequel of the paper, this common object wil be referred to as $(\Et_{\bullet}, \tt )$. 
\end{defn}

\begin{proof}
It follows from the results discussed in Theorem \ref{openint} that on the open set $\Ct \setminus \Pt$ 
the two Higgs bundles $(\Et^{\prime},\tt^{\prime})$ and $(\Et,\tt)$ agree. Indeed, over $\Ct \setminus \Pt$ the 
two surfaces $\CP \times \CPt$ and $\Zt^\Pt$ are isomorphic, the same holds for the sheaves $M^{0}$ and
$\Mt^\Pt$ and the coordinates $\xi$ and $\Xt_\Pt$. Finally, the two factors of $1/2$ in the 
definition of $(\Et,\tt)$ (namely, that of $\theta_{\xi}=\theta-\xi/2$ and $\tt=-z(\xi)/2$) cancel each other. 
Hence, we only need to check that the extensions to the singularities agree as well. 

Lemma \ref{lem:Purity} implies
$$
  (\pmapz)_{\ast} \circ (\qmapP)_*(N^{\abs})=N^{0}(\pi^{-1}(\infty))
$$ 
as sheaves, and because 
$\widehat{\sE}^{-}\cap\Sigma^{\abs}=(\pi \circ 
\pmapz\circ \qmapP)^{-1}(\infty)\cap\Sigma^{\abs}$ and 
$(\pit \circ \pmapz\circ \qmapP)^{-1}(\Pt)=
\widehat{\sE}^{-}\cup\widehat{\sE}^{+}$, this implies 
$$
    (\pmapz)_{\ast} \circ (\qmapP)_*\add{\widehat{\sE}^{+}}(N^{\abs})=N^{0}(\pit^{-1}(\Pt)). 
$$
Since the projections $\pit \circ \pmapz \circ \qmapP$ and 
$\pit \circ \pmaptP \circ \qmapPt$ from $Z^{\abs}$ to $\CPt$ are the same, we have the isomorphism of sheaves 
$$
     \pit_{\ast} N^{0}(\Pt)= \pit_{\ast} \Mt^\Pt, 
$$
because both are equal to the direct image of $\add{\widehat{\sE}^{+}}N^{\abs}$ with respect to the same projection. 
However, the direct image of the parabolic structure of $N^{0}$ is not the correct one: indeed, 
on one hand the set of parabolic points of $N^{0}$ contains the points of $\sE^{+}\cap \Sigma^{0}$ with trivial 
parabolic structure, so these will induce extra parabolic points on $\CPt$ with trivial structure; 
and on the other hand it does not contain the points $\widehat{\sE}^{+}\cap \Sigma^{0} \subset \pit^{-1}\Pt$ so that the direct 
image with respect to $\pit$ of the parabolic structure on $N^{0}$ does not really make sense. 
On the other hand, we modified the parabolic divisor of $N^{\abs}$ so that these problems do not 
occur when we push it down. 
Hence, in order to prove equality of the bundles $\Et^{tr}$ and $\Et^{\prime}$ it is sufficient to prove that 
$\pi_{\ast}N^{0}=\Et^{tr}$; whereas for their parabolic structure, we will to work directly with 
$\Mt^\Pt(-\pit^{-1}(\Pt))$. 

Now, as local sections of $N^{0}$ are sections of $M^{0}$ vanishing in the points 
$T^{-}=(\{\ti\})\cap \Sigma^{0}$ and $\Tt^{-}=(\Pt \times\{\infty \})\cap \Sigma^{0}$, this means precisely that if 
$\zeta$ is a local coordinate of 
$\CPt$ at $\ti$ then the local sections near $\ti$ of the direct image $\pit_{\ast}N^{0}$ can be represented by 
a local section of the sheaf $M^{0}$ multiplied by bump-functions concentrated at the spectral points
of $\zeta$ whose heights converge to $0$ up to first order as $\zeta \ra 0$. On the other hand, the induced extension
at infinity is defined precisely by admitting a representation by bump-functions of constant height as 
$\zeta \ra 0$, and the local sections of the transformed extension are obtained from those of the induced
extension upon multiplication of these latter with $\xi^{-1}=\zeta$ (c.f. the discussion before \ref{degformula}). 
It is proved in Proposition 4.24 of \cite{Sz} that a $D''_{\xi}$-harmonic $1$-form 
$\varphi =\varphi_{1}\d z + \varphi_{\bar{1}}\d \bar{z}$ represents the element of
$M^{0}_{\xi}$ which is the class of $\{\varphi_{1}(q(\xi))\d z\}$ modulo the image of $\theta$, 
where $q(\xi)$ runs over the finite set of spectral points of $\xi$. It follows that 
multiplying the harmonic representatives 
by bump-functions of height converging to $0$ instead of constant height amounts to taking sections of 
$M^{0}$ that vanish at $\ti$. 
Therefore, locally near the dual infinity the isomorphism of holomorphic bundles $\pit_{\ast}N^{0}=\Et^{tr}$ holds. 
Similarly, near a logarithmic singularity $\xi_{l}$ the change of trivializations to obtain $\Et^{tr}$ from 
$\Et^{ind}$ is to take bump-functions of height decaying as $|\xi -\xi_{l}|$ near the spectral points converging
to $\infty \in \CP$. This amounts 
precisely to taking local sections of $M^{0}$ vanishing up to first order on the divisor $\{\infty \} \times \Pt$. 
Such local sections of $M^{0}$ are by definition the local sections of $N^{0}$, 
therefore the direct image of $N^{0}$ in the logarithmic singularities $\Pt$ of $\tt$ is also equal to
$\Et^{tr}$; this implies isomorphism of the bundles and Higgs fields. 

It remains to identify the parabolic structures of the direct image. First of all, the set of parabolic points
of $\Et^{tr}_{\bullet}$ in $\CPt$ is $\Dt=\Pt \cup \{\ti \}$, and the same holds for $\Et^{\prime}_{\bullet}$ 
because the deletion procedure removes extra parabolic points with trivial structure from $\Mt^\Pt_{\bullet}$. 
Second, by Section 4.6 of \cite{Sz}, near the 
punctures the local bases of the transformed bundle defined by the representatives given in the previous
paragraph are adapted to the transformed harmonic metric. Hence, the direct images of the parabolic filtrations
of $\Mt^\Pt(-\pit^{-1}(\Pt))_{\bullet}$ are the filtrations of $\Et^{tr}_{\bullet}$. Furthermore, the same thing holds 
for the weights as well. Indeed, the weights of $\Mt^\Pt(-\pit^{-1}(\Pt))_{\bullet}$ in the points of 
$\widehat{\sE}^{-}\cap\widehat{\Sigma}^{\Pt}$ above $\xi_{l}$ are equal to the parabolic weights $\alpha^{\infty}_{k}$ of the
original bundle on the $\xi_{l}$-eigenspace of $A$ at infinity. This follows because the weights of $M^\sP_{\bullet}$ at 
infinity are $\alpha^{\infty}_{k}$, because near infinity the isomorphism of sheaves $N^\sP\cong M^\sP$ holds, and 
because our convention is to keep the same parabolic weights for sheaves isomorphic near a parabolic divisor. 
The weights of $\Mt^\Pt(-\pit^{-1}(\Pt))_{\bullet}$ in the points $\Tt^{+}$ are in turn equal to $0$ by the
definition of adding a new divisor to the parabolic divisor of a parabolic structure. 
On the other hand, by Theorem 4.37 of \cite{Sz} the non-zero weights of $\Et^{tr}_{\bullet}$ in $\Pt$ corresponding 
to the sections defined above are equal to $\alpha^{\infty}_{k}$. This proves equality of the parabolic structures 
of the two extensions in the logarithmic singularities. 
Similarly, the weights of $\Mt^\Pt(-\pit^{-1}(\Pt))_{\bullet}$ at $\ti$ are the non-zero weights $\alpha^{j}_{k}$ of 
$\E_{\bullet}$ in points of $\sP$: again, this follows from the fact that the weights of $M^\sP_{\bullet}$ in 
the points of 
$\pi^{-1}(\sP)\cap \Sigma^\sP\setminus T^{+}$ are $\alpha^{j}_{k}$, combined with the local isomorphism of sheaves 
$N^\sP\cong M^\sP$ away from the $0$-section of $\pi$. By Theorem 4.34 of \cite{Sz} the corresponding weights 
of $\Et^{tr}_{\bullet}$ in $\ti$ are also equal to $\alpha^{j}_{k}$; whence the theorem. 
\end{proof}

\section{Examples}\label{sec:Ex}
In this section, we illustrate by two examples how our approach allows to increase the degree of generality 
of the setup of the transform defined in \cite{Sz}. 

\subsection{Nilpotent residues}\label{nilpsect}
In this first example we show that the transform can be defined for a Higgs field whose polar parts 
are not necessarily semi-simple. The conclusion is that a zero residue matrix at infinity 
can induce two nilpotent residues of rank one in points of finite distance of the transformed object -- 
hence, there is no analog to nilpotent parts of the preservation of the sum of the ranks of the semi-simple 
parts of the residues by the transform. Paralelly to this, the multiplicity of the parabolic weights 
is also not preserved by the transform. In concrete terms, this means that a parabolic weight $\alpha$ of 
multiplicity $1$ splits up to a multiplicity $2$ weight $\alpha/2$; in particular, the total parabolic degree 
is preserved. Notice finally that in this example we start out with a Higgs field with a rank-two 
singularity at infinity, and we arrive at one with a logarithmic pole at infinity. 

Let $u_{0}$ and $v_{0}$ be the standard coordinates on $\CP$ in a neighborhood of $0$ and $\infty$ respectively, 
$\sP=\{0\}\subset \CP$, $\E$ be the rank-two trivial holomorphic bundle $\O_{\CP}\oplus\O_{\CP}$, and $\theta$ on the open
affine $v_{0}=1$ containing $0$ is given in matrix form 
\begin{eqnarray}\label{nilpfield}
     \left (
     \begin{array}{rr}
            \frac{1}{u_{0}} & 1 \\
            -1 & -\frac{1}{u_{0}}
     \end{array}\right );
\end{eqnarray}
or, in homogeneous coordinates, 
\begin{eqnarray}\label{nilpfield2}
     \left (
     \begin{array}{rr}
            v_{0} & u_{0} \\
            -u_{0} & -v_{0}
     \end{array}
     \right ).
\end{eqnarray}
The residue in $0$ has two distinct eigenvalues $\pm 1$; whereas the limit of the field at infinity is 
\begin{eqnarray*}
     \left( 
     \begin{array}{rr}
            0 & 1 \\
            -1 & 0
     \end{array}
     \right ),
\end{eqnarray*}
with eigenvalues $\pm i$. Furthermore, setting $u_{0}=1$ in (\ref{nilpfield2}), an easy computation shows that 
the eigenvalues $\xi(v_{0})$ of the matrix can be written 
\begin{equation}\label{eigenvalues}
   \xi(v_{0})=\pm i\sqrt{1-{v_{0}^{2}}}=\pm i \left( 1- \frac{v_{0}^{2}}{2}+ O(v_{0}^{4}) \right).
\end{equation}
In particular, since the eigenvalues are distinct for $v_{0}=0$, in a neighborhood of $\infty$ there exists a
trivialization of $\E$ in which the matrix of this endomorphism is diagonal. 
This trivialization then clearly satisfies the properties required by (\ref{thetainf})-(\ref{binf}) with 
first-order term $B_{\infty}=0$, since the eigenvalues are functions of ${v_{0}^{2}}$. However, 
the assumption that the eigenvalues of $B_{\infty}$ are all non-vanishing and distinct obviously fails. 
Finally, let $\alpha^{0}_{+},\alpha^{0}_{-}\in [0,1[$ be arbitrary 
weights at the singularity $0$ corresponding to the $1$ and $(-1)$-eigenspaces respectively; 
and let $\alpha^{\infty}_{+},\alpha^{\infty}_{-}\in [0,1[$ be arbitrary weights at infinity corresponding to the $i$ and 
$(-i)$-eigenspaces. 

Consider the standard spectral surface $Z^{1}=\PP_{\CP}(\O_{\CP}\oplus\O_{\CP}(\sP))$, and call 
 the preferred sections of  $\O_{Z^{1}}(1)\otimes\O_{\CP}(\sP)$ and $\O_{Z^{1}}(1)$, $x_{1},y_{1}$ respectively. 
The standard spectral curve $\Sigma^{1}$ in $Z^{1}$ is defined by the polynomial 
$$
     \det(x_{1}-y_{1}\theta).
$$
Since the standard spectral curve does not intersect the infinity-section, we may assume $y_{1}=1$. 
Then this polynomial becomes 
$$
     x_{1}^{2}-v_{0}^{2}+u_{0}^{2}. 
$$
The solution of this homogeneous equation is 
\begin{eqnarray}\label{emb1}
    \begin{array}{cccc}
    y_{1}=1 & x_{1}=s^{2}-t^{2} & u_{0}=2st & v_{0}=s^{2}+t^{2}, 
    \end{array}
\end{eqnarray}
where $[s:t]$ stands for homogeneous coordinates on a smooth rational curve. In other words, the mapping  $[s:t]\mapsto u_0,v_0,x_1,y_1$ defined by (\ref{emb1})
is a closed embedding whose image is the smooth spectral curve $\Sigma^{1}$. 
In particular, it has two branches over $[u_{0}=0:v_{0}=1]=0\in \C$ which do not intersect: 
one of them through $s=0,t=1$, the other one through $s=1,t=0$. The first corresponds to $x_{1}=-1$, 
that is the eigenvalue $-1$ of the residue of $\theta$ in $0$, the second $x_{1}=1$ to the eigenvalue $1$. 
Similarly, the two branches of $\Sigma^{1}$ over $\infty\in \CP$ pass through $x_{1}=i$ and $x_{1}=-i$. 
The pull-back of $\E$ to $\Sigma^{1}$ is the rank-two trivial holomorphic bundle $\O_{\Sigma^{1}}\oplus\O_{\Sigma^{1}}$
and the map $\Theta$ is then by definition 
$$
   x_{1}-y_{1}\theta: \O_{\Sigma^{1}}\oplus\O_{\Sigma^{1}} \lra \O_{\Sigma^{1}}(2)\oplus\O_{\Sigma^{1}}(2). 
$$
Using the above formulae, we obtain for this map the matrix form 
$$
    \left (
    \begin{array}{rr}
        -2t^{2} & -2st \\
        2st & 2s^{2}
    \end{array} \right ).
$$
A cokernel map for this is left matrix multiplication 
$$
   (s,t): \O_{\Sigma^{1}}(2)\oplus\O_{\Sigma^{1}}(2)\lra \O_{\Sigma^{1}}(3);
$$
in particular, the sheaf $M^\sP$ is $\O_{\Sigma^{1}}(3)$. 
Now, since the residue of $\theta$ in $0$ has two distinct non-zero eigenvalues $\{\pm 1\}$ and
the rank of $\E$ is equal to $2$, the set $t^{=0}$ is empty. Therefore, the sheaf $N^\sP$
is the kernel of evaluation of $\O_{\Sigma^{1}}(3)$ in the points of the spectral curve $\Sigma^{1}$ over
the point at infinity $[u=0:v=1]$. Because $\Sigma^{1}$ is a double cover of $\CP$ and smooth over infinity, 
we deduce that $N^\sP=\O_{\Sigma^{1}}(1)$. Furthermore, it has four parabolic points: the two points over 
$0\in \C$ and the two points over $\infty \in \CP$ discussed above. The weights are as follows: the one in $x_{1}=-1$ 
over $0 \in \C$ is $\alpha^{0}_{-}$; the one in $x_{1}=1$ over $0 \in \C$ is $\alpha^{0}_{+}$; the one in $x_{1}=-i$ 
over $\infty \in \CP$ is $\alpha^{\infty}_{-}$; the one in $x_{1}=i$ over $\infty \in \CP$ is $\alpha^{\infty}_{+}$. 

Let us now consider the compactification $\Sigma^{0}$ of the open spectral curve 
$\Sigma^\oo$ in the surface $\CP \times \CPt$. 
We have seen that it is the proper transform of $\Sigma^{1}$ with respect to the elementary transformation 
linking $Z^{1}$ to $\CP \times \CPt$. Let $x_{0},y_{0}$ be homogeneous coordinates of $\CPt$: they can be thought
of as sections of $\O_{\CPt}(1)$ vanishing in $0$ and $\ti$ respectively. Since the elementary transformation 
in question blows up the point $u_{0}=0$, the relation between these coordinates and the sections of $\O_{Z^{1}}(1)$ is 
\begin{align}
     x_{0}=x_{1} && y_{0}=y_{1}u_{0}. 
\end{align}
Therefore, the parametrization (\ref{emb1}) transforms into 
\begin{eqnarray}\label{emb0}
    \begin{array}{cccc}
    y_{0}=2st & x_{0}=s^{2}-t^{2} & u_{0}=2st & v_{0}=s^{2}+t^{2},
    \end{array}
\end{eqnarray}
and the equation defining the curve becomes 
\begin{eqnarray}\label{spcurveqn0}
     x_{0}^{2}u_{0}^{2}-y_{0}^{2}v_{0}^{2}+y_{0}^{2}u_{0}^{2}.
\end{eqnarray}
This curve in $\CP \times \CPt$ is not smooth. Indeed, it is straightforward to check that 
it has a node in the point $(0,\ti)\in \CP \times \CPt$. On the other hand, it has no other singularity, 
because on the complementary of the fiber of $\pi$ over $0$ it is isomorphic to $\Sigma^{1}$. 

The first thing we need to identify in order to perform the elementary transformations for the 
projection $\pit$ is the polar divisor $\Pt\subset \CPt$. 
Recall that it is given as the 
intersection points of $\Sigma^{1}$ with the $\infty$-fiber of $\pi$. Plugging $u_{0}=0,v_{0}=1$ in the 
equation of $\Sigma^{1}$ we get $x_{1}^{2}=-1$. We deduce $\Pt=\{[i:1],[-i:1]\}$, in other words 
the points $\{\pm i\}$ of $\Ct$. We now proceed in two steps: first, compute the coordinates of 
the surface $\Zt^{i}$ obtained by performing an elementary transformation in the point $[i:1]$; 
then perform another elementary transformation, this time in $[-i:1]$, to obtain the surface $\Zt^2=\Zt^{\{ -i,i\}}$. 

The transformation in $i$ is 
\begin{align*}
   u_{1}=u_{0}(x_{0}-iy_{0})=u_{0}(s-it)^{2} && v_{1}=v_{0};
\end{align*}
since these are projective variables, this is equivalent to 
\begin{align*}
   u_{1}=u_{0}(s-it) && v_{1}=\frac{v_{0}}{s-it}=s+it. 
\end{align*}
The transformation in $-i$ is 
\begin{align*}
   u_{2}=u_{1}(x_{0}+iy_{0})=u_{1}(s+it)^{2} && v_{2}=v_{1}, 
\end{align*}
which is again equivalent to 
\begin{align}\label{emb2}
   u_{2}=u_{1}(s+it)=u_{0}(s-it)(s+it)=2st(s^{2}+t^{2}) && v_{2}=\frac{v_{1}}{s+it}=1.
\end{align}
Meanwhile, since over the points $\pm i \in \Ct$ the spectral curve has only one branch through $\infty \in \CP$, 
these modifications do not introduce new singularities. In other words, the spectral curve
is transformed by these modifications into a rational curve $\widehat{\Sigma}^2=\widehat{\Sigma}^{-i,i}$ 
with one node, in the zero section over the point $\ti$. 

Next, we need to compute the transformed bundle $\Et$. First, because $\Sigma^1$ does 
not pass through the intersection of the $0$-section and the $0$-fiber of $Z^1$, the curve $\Sigma \subset Z$ 
is isomorphic to $\Sigma^1$, and $\Sigma$ is disjoint from the proper transform $\sE^+$ of the blow-up $Z \to Z^1$. 
Furthermore, because the multiplicity of each intersection point of $\Sigma^1$ 
and the infinity-fiber of $Z^1$ is $1$, the curve $\Sigma^{\abs} \subset Z^{\abs}$ is also isomorphic to $\Sigma^1$, which 
is a projective line. It follows that the proper transform of $M^1$ in $Z^{\abs}$ is $\O_{\Sigma^{\abs}}(3)$, 
and the intermediate spectral sheaf is by definition $N^{\abs}=\del{\sE^+}(\O_{\Sigma^{\abs}}(3))$, which is 
$\O_{\Sigma^{\abs}}(3)$ because $\sE^+$ is disjoint from $\Sigma^{\abs}$. The curve $\widehat{\Sigma}^2$ intersects 
the zero-section of $\Zt^2$ over each of the points $\pm i$ with multiplicity $1$. Hence, both components of the 
exceptional divisor $\widehat{\sE}^+$ of the blow-up of $\Zt^2$ in these points intersects the curve $\Sigma^{\abs}$ 
in one point. Therefore, the sheaf $\add{\widehat{\sE}^+}(N^{\abs})=(\O_{\Sigma^{\abs}}(3))(\widehat{\sE}^+)$ 
is isomorphic to $\O_{\Sigma^{\abs}}(5)$. By definition, $\Et(\{ -i,i\})$ is the direct image of 
$\add{\widehat{\sE}^+}(N^{\abs})$ with respect to the projection of $Z^{\abs}$ to $\CPt$. By the computations above 
and because $\Sigma^{\abs}$ is a double cover of $\CPt$, this means that $\Et$ is the direct image of 
$\O_{\Sigma^{\abs}}(1)$. We claim that 
$$
    \Et=\O_{\CPt}\oplus\O_{\CPt}.
$$
Indeed, since $\Sigma^{\abs}$ is a double cover of $\CPt$, clearly $\Et$ is of rank $2$. 
Now $\O_{\Sigma^{\abs}}(1)$ has exactly two independent global sections: $s$ and $t$. They induce 
global sections of $\Et$. Conversely, any global section of $\Et$ induces a global section 
of $\O_{\Sigma^{\abs}}(1)$, and is therefore a linear combination of $s$ and $t$. In different 
terms, $s$ and $t$ give a global trivialization of $\Et$ on $\CPt$. 

The last thing to compute is the transformed Higgs field $\tt$. Here, 
the section $\Xt_{\Pt}$ is called 
 $u_{2}$, and $\tt$ is the direct image of multiplication by $-u_{2}$ on $\O_{\widehat{\Sigma}^{-i,i}}(1)$. 
Notice that using (\ref{emb0}) and (\ref{emb2}) we obtain
\begin{align}
    u_{2}\cdot s&=2st(s^{2}+t^{2})\cdot s=x_{0}y_{0}\cdot s + y_{0}^{2}\cdot t\\
    u_{2}\cdot t&=2st(s^{2}+t^{2})\cdot t=y_{0}^{2}\cdot s - x_{0}y_{0}\cdot t, 
\end{align}
therefore the matrix of $\tt$ in the above trivialization is 
$$
    -\left( \begin{array}{rr}
         x_{0}y_{0}&  y_{0}^{2}\\
        y_{0}^{2} & - x_{0}y_{0}
    \end{array} \right ).
$$
Here $x_{0},y_{0}$ are standard projective coordinates for $\CPt$, vanishing at the points $0$ and 
$\ti$ respectively. The matrix form of this map on the affine $\Ct$ with poles in the points 
$\{ \pm i\}$ can be obtained by setting $y_{0}=1$, and dividing each entry by $(x_{0}-i)(x_{0}+i)$. 
The result is 
$$
   \tt = -\frac{1}{x_{0}^{2}+1}\begin{pmatrix}
         x_{0}&  1\\
        1 & - x_{0}
    \end{pmatrix}.
$$
This matrix clearly has poles in $x_{0}=\pm i$, with residues 
$$
    -\begin{pmatrix}
        \pm i &  1\\
        1 & \mp i 
    \end{pmatrix},
$$
both of which are nilpotent. Furthermore, the spectral curve is ramified over both of these points, as it 
can be seen for example from the defining equation (\ref{spcurveqn0}) of $\Sigma^{0}$ upon putting $x_{0}=1$, 
$v_{0}=1$ and using the observation that taking proper transform with respect to an elementary transformation 
does not change the property of the projection to $\CPt$ being ramified or not. Another way of seeing the same
thing is as follows. Express $v_{0}$ in terms of $x_{0}=\xi$ in (\ref{eigenvalues}) for example near the 
value $x_{0}=i$; we obtain that 
\begin{equation}\label{vform}
    v_{0}^{2}=x_{0} - i, 
\end{equation}
which is clearly an index $2$ ramification for $\pit$. It follows that over the logarithmic poles $\pm i$ 
of the transformed field, both branches of the spectral curve pass through the infinity-section of $\pit$. 

We deduce that the parabolic filtration of $\Et$ in these points has to be trivial, hence with only one weight
in $i$ (resp. $-i$). Moreover, the norm squared with respect to $h$ of the cokernel vector is equivalent to 
$|v_{0}|^{2\alpha^{\infty}_{+}}$, which is equal to $|x_{0}-i|^{\alpha^{\infty}_{+}}$ because of (\ref{vform}); therefore, 
this unique parabolic weight in the point $i$ (resp. $-i$) is $\alpha^{\infty}_{+}/2$ (resp. $\alpha^{\infty}_{-}/2$). 
On the other hand, near $\ti$ the matrix for $\tt$ looks up to higher-order terms like 
$$
    -\begin{pmatrix}
        \frac{1}{x_{0}} &  \frac{1}{x_{0}^{2}}\\
        \frac{1}{x_{0}^{2}} & - \frac{1}{x_{0}}
    \end{pmatrix}; 
$$
this converges to $0$ as $x_{0}$ goes to infinity, and the first-order term in its Taylor series is 
$$
    \begin{pmatrix}
        - 1 &  0\\
        0 & 1
    \end{pmatrix},
$$
with eigenvalues $\{\pm 1\}$. The parabolic weight of the $\pm 1$-eigenspace is $\alpha^{0}_{\pm}$.

\subsection{Higher order pole}\label{hopsect}
Although so far we assumed that the Higgs field has at most logarithmic singularities in the 
singular points at finite distance, it is relatively clear that iterating the construction several 
times according to the order of the poles, one can get a transform for Higgs bundles with higher-order poles -- 
the transformed Higgs field will then have a ramification at infinity. 
Here we describe the archetype of this phenomenon: the original Higgs bundle has a maximal 
ramification at infinity, and the transformed field has a pole in the origin whose order equals 
to the index of this ramification. 

Let $r\geq2$ and take $\E$ to be the rank $r$ trivial holomorphic bundle $\O_{\CP}^{\oplus r}$, with $\sP=\{0\}$ the only 
regular parabolic point. Define the Higgs field as the map 
\begin{align*}
     \theta:\O_{\CP}^{\oplus r}\lra \O_{\CP}(1)^{\oplus r} 
\end{align*}
defined in some global trivialization $\{\tau_{1},\ldots ,\tau_{r}\}$ by the matrix 
\begin{align}\label{mat:hop}
     \begin{pmatrix}
          0 & u_{0}&0&\ldots &0\\
          0 & 0 & u_{0}&\ldots &0 \\
          \vdots  & \vdots &  \ddots & \ddots & \vdots  \\
          0 & 0  & \ldots   & 0 & u_{0}\\
          v_{0} & 0 & \ldots & 0 & 0
     \end{pmatrix},
\end{align}
where $u_{0},v_{0}$ are the global sections of $\O_{\CP}(1)$ vanishing in $0$ and $\infty$ respectively. 
Here and in all this section we identify $\O_{\CP}(\sP)$ with $\O_{\CP}(1)$, and correspondingly replace $\sP$ by
$1$ in all upper and lower indices. 

We immediately see that at infinity the matrix (\ref{mat:hop}) has only $0$ eigenvalues, so we deduce that the singular set
of the transform will be $\Pt=\{0\}$. 
Furthermore, in a local affine coordinate $v$ centered at infinity, the matrix becomes 
\begin{align*}
     \begin{pmatrix}
          0 & 1&0&\ldots &0\\
          0 & 0 & 1&\ldots &0 \\
          \vdots  & \vdots &  \ddots & \ddots & \vdots  \\
          0 & 0  & \ldots   & 0 & 1\\
          v & 0 & \ldots & 0 & 0
     \end{pmatrix}.
\end{align*}
It is clear that there exists no non-trivial subspace invariant by both the constant and the first-order term of this matrix. 
Since the parabolic filtration has to be preserved by the polar part of the field, this then implies that the only possible 
filtration in this point is the trivial filtration 
$$
   \E|_{\infty}=F_{0}\E|_{\infty}\supset F_{1}\E|_{\infty}=\{0\}.
$$ 
Let us call the corresponding weight $\alpha^{\infty}$. By the general hypotheses made on the weights, $\alpha^{\infty}$ 
is in $]0,1[$. For the sake of simplicity, let us also suppose that $\alpha^{\infty}<1/r$. 

On the other hand, the residue of the Higgs field (\ref{mat:hop}) in the only regular 
singular point $0\in\C$ is of rank $1$, so the transformed bundle will be of rank $1$. Moreover, since this 
residue has only $0$ eigenvalues, the standard spectral curve passes through the zero-section over the 
polar point $0$ with maximal multiplicity $r$. Condition \ref{assn:Main} then forces that the parabolic 
structure in this polar point is trivial, i.e. the filtration is trivial and the only weight is $0$. 

The spectral surface is $Z^{1}=\PP_{\CP}(\O_{\CP}\oplus\O_{\CP}(1))$, and 
\begin{align*}
   \Theta_{1}= x_{1}-y_{1}\theta = \begin{pmatrix}
      x_{1} & -y_{1}u_{0}&0&\ldots &0\\
          0 & x_{1} &  -y_{1}u_{0}&\ldots &0 \\
          \vdots  & \vdots &  \ddots & \ddots & \vdots  \\
          0 & 0  & \ldots   & x_{1} &  -y_{1}u_{0}\\
           -y_{1}v_{0} & 0 & \ldots & 0 & x_{1}
      \end{pmatrix}
\end{align*}
implies that the spectral curve is 
$$
   \Sigma^{1}=(\mbox{det}(\Theta_{1}))=(x_{1}^{r}-(-y_{1})^{r}u_{0}^{r-1}v_{0}).
$$
This curve is singular in the point $x_{1}=0,u_{0}=0$ if $r>2$. 
Therefore, instead of working on this surface, we choose first to perform an elementary transformation 
in the point $0$ and reduce our problem to the case of a smooth curve in $\CP \times \CPt$. 
The elementary transformation to apply is given by the coordinate changes $x_{0}u_{0}=x_{1}, y_{0}=y_{1}$, 
and now we consider $Q^{-1}\theta$ as a map $\E \ra \F(1)=\O_{\CP}^{\oplus(r-1)}\oplus\O_{\CP}(1)$. 
In concrete terms, denoting by $\boxtimes$ external tensor product of sheaves on a product space, 
the map $Q$ in the Diagram (\ref{diag:Naive3}) has the form $\mbox{diag}(u_{0},\ldots ,u_{0},1)$ in the same basis as above, 
and this means that 
$$
   \Theta_{0}=Q^{-1}(x_{0}u_{0}-y_{0}\theta):\O_{\CP}^{\oplus r}\lra(\O_{\CP}^{\oplus(r-1)}\oplus\O_{\CP}(1))\boxtimes\O_{\CPt}(1)
$$
is given by 
\begin{align*}
   \begin{pmatrix}
      x_{0} & -y_{0}&0&\ldots &0\\
          0 & x_{0} &  -y_{0}&\ldots &0 \\
          \vdots  & \vdots &  \ddots & \ddots & \vdots  \\
          0 & 0  & \ldots   & x_{0} &  -y_{0}\\
           -y_{0}v_{0} & 0 & \ldots & 0 & x_{0}u_{0}
      \end{pmatrix}.
\end{align*}
Therefore, the spectral curve $\Sigma^{0}$ is defined by the equation 
$$
     x_{0}^{r}u_{0}-(-y_{0})^{r}v_{0}=0, 
$$
and this is clearly a non-singular rational curve smoothly parametrized by $(x_{0},y_{0})$: namely, 
one has $u_{0}=(-y_{0})^{r}, v_{0}=x_{0}^{r}$. 
As an effect of passing to the product surface we therefore desingularize the curve, 
and applying the map $Q^{-1}$ we get rid of the extra fiber of multiplicity $(r-1)$ over $0$ of the total 
transform of the curve. 
Over the origin in $\CP$ one has $u_{0}=0,v_{0}=1$, so necessarily $y_{0}=0$, which means that the only point 
of the spectral curve over the origin is the point $(0,\ti)\in \CP \times \CPt$. Moreover, putting $x_{0}=1$ the 
equation for the curve near this point becomes 
$$
   u_{0}=(-1)^{r}y_{0}^{r},
$$
which means that the projection to $\CP$ has a ramification of index $r$ in this point. 
Similarly, the curve passes through $(\infty,0)$ and projection to $\CP$ has a ramification of index $r$ 
over the infinity. In particular, it intersects the $0$- and $\infty$-fibers of $\pi$ in these points with 
multiplicity $r$. A simple computation shows that the map 
\begin{align*} 
   A&:(\O_{\CP}^{\oplus(r-1)}\oplus\O_{\CP}(1))\boxtimes\O_{\CPt}(1) \lra \O_{\CP}(1)\boxtimes \O_{\CPt}(r)\\ 
   A&=(y_{0}x_{0}^{r-2}v_{0},y_{0}^{2}x_{0}^{r-3}v_{0},\ldots, y_{0}^{r-1}v_{0},x_{0}^{r-1})
\end{align*}
is a cokernel for $\Theta_{0}$. In particular, the cokernel sheaf $M^{0}$ of $\Theta_{0}$ is the restriction of 
$\O_{\CP}(1)\boxtimes \O_{\CPt}(r)$ to $\Sigma^{0}$, which is equal to $\O_{\Sigma^{0}}(2r)$ because $\Sigma^{0}$ 
is an $r$-to-$1$ cover of $\CP$ and a $1$-to-$1$ cover of $\CPt$. By definition, the sheaf $N^{0}$ 
is the kernel of evaluation of $M^{0}$ in the intersection points  $(\infty,0)\in \CP\times \CPt$ and 
$(0,\ti)\in \CP\times \CPt$ of the spectral curve with the $\infty$-fiber of $\pi$ and the $0$-fiber in the 
$\infty$-section. We have seen above that these intersections are of multiplicity $r$. It follows that 
$N^{0}=M^{0}(-2r)=\O_{\Sigma^{0}}$. Since $\pit:\Sigma^{0}\ra\CPt$ is an isomorphism, we deduce that $\Et=\O_{\CPt}$. 

Let us now identify the transformed Higgs field $\tt$; for this purpose, we need to perform additional 
elementary transformations on $\CP \times \CPt$, but this time with respect to the projection $\pit$. 
Namely, since the spectral curve $\Sigma^{0}$ intersects the $\infty$-section of $\pit$ in its $0$-fiber, we need 
to introduce $u_{1}=u_{0}x_{0},v_{1}=v_{0}$. The equation of the proper transformed curve $\widehat{\Sigma}^1$ 
is then given by $x_{0}^{r-1}u_{1}-(-y_{0})^{r}v_{1}=0$. However, this still intersects the $\infty$-section of 
$\pit$, so we need to do another elementary transformation $u_{2}=u_{1}x_{0},v_{2}=v_{1}$, and continue this 
procedure until the proper transformed curve no more intersects the $\infty$-section, that is 
$u_{r}=u_{0}x_{0}^{r},v_{r}=v_{0}$. The equation of the proper transformed curve $\widehat{\Sigma}^r$ is now 
$u_{r}-(-y_{0})^{r}v_{r}=0$. Since this curve does not intersect the $\infty$-section of $\pit$, we may set
$v_{r}=1$. Then the curve is given by $u_{r}=(-y_{0})^{r}$. Now, one has by definition 
$\tt=\pit_{\ast}(-u_{r}\cdot)$, hence we see that the transformed Higgs field has the form 
$$
   -(-y_{0})^{r}:\O_{\CPt}\lra\O_{\CPt}(r),
$$
where $\O_{\CPt}(r)$ stands for $\O_{\CPt}(r\{0\})$. This map therefore has an order $r$ pole in $0$ 
(on the affine $\Ct$ it can be written as $\pm1/x_{0}^{r}$), and on the other hand it clearly has an order 
$r$ zero at infinity. The fibers being of dimension one, the parabolic filtrations are trivial in both of 
these points. Similar arguments as in Subsection \ref{nilpsect} show that the corresponding parabolic 
weights in $0$ and $\ti$ are $r\alpha^\infty$ and $0$ respectively.

\section{Quasi-involutibility}\label{sec:qi}

As the second author has proved it in Chapter 5 of \cite{Sz}, one of the features of Nahm transform is its 
involutibility up to a sign. The proof there is done in the framework of integrable connections, and 
relies on the analysis of a spectral sequence. Our aim in this section is to 
give a new, more geometric proof of the same result in terms of the techniques developped in this paper
(Theorem \ref{invthm}). 

Define the $\kk$-linear map 
$$
    -1:\CP \lra \CP
$$
to be the extension to $\CP$ of the ($\kk$-linear) involution taking an element of $\kk$ to its additive
inverse. It has two fixed points: $0$ and $\infty$. We will denote by $(-1)_{rel}$ the relative version of this map
on any $\CP$-fibration over a curve. It then induces a map on any blow-up of points in the $0$-section or the 
$\infty$-section, that we will still denote with the same symbol. 
\begin{rk}
In what follows, it will be important to distinguish the divisors $(-1)^{\ast}\sP$ and $-\sP$: the first is the set of points 
$-p$ where $p\in \sP$, whereas the second is the inverse of the divisor $\sP$ in the divisor group. 
\end{rk}

Recall that given a parabolic Higgs bundle $(\E,\theta)$ on $\CP$ with singularities on $D \cup \{\infty\}$ we have constructed 
in Section \ref{trsect} its Nahm transformed Higgs bundle $(\Et,\tt)$: it is a parabolic Higgs bundle on $\CPt$, 
the "dual" projective line, with singularities on $\Dt =\Pt\cup\{\ti\}$. Here $\Pt$ is the set of eigenvalues of the leading
term of $\theta$ at $\infty$. We have also computed the eiganvalues of $\tt$ at $\ti$, and we realized that they 
agree with the image $-1(\sP)$ of $\sP$ under the involution. 
The bidual $\CPtt$ of $\CP$ identifies naturally with $\CP$ itself, hence applying the Nahm transform to 
$(\Et,\tt)$, we obtain a parabolic Higgs bundle $(\Ett,\ttt)$ on $\CP$ with singularities on the set $-1(\sP)\cup \{\infty\}$. 
The main result of this section can now be formulated. 
\begin{thm}\label{invthm}
If $(\E,\theta)$ satisfies  Condition \ref{assn:Main}, then 
there is a natural isomorphism of parabolic Higgs bundles between $\left(\Ett,\ttt \right)$ and $(-1)^{\ast}(\E,-\theta)$. 
\end{thm}
\begin{rk}
There is a sign change of $\theta$ between this and the corresponding formula in \cite{Sz}. This is because 
there we considered the Higgs field as a $1$-form valued endomorphism, and $\d(-z) =-\d z$. 
\end{rk}
\begin{proof}

Starting with the parabolic Higgs bundle $(\Et,\tt)$ on $\CPt$, we wish to construct its transform. The first
object we need to understand is the standard spectral triple $(W^{\Pt},\Xi^{\Pt},Q^{\Pt})$ of $(\Et,\tt)$. Remember
that in Section \ref{trsect} we constructed the spectral triple $(\Zt^{\Pt},\widehat{\Sigma}^{\Pt},\Mt^{\Pt}_{\ast})$
out of $(\E,\theta)$. First, because of the definition 
$W^{\Pt}={\mathbf P}(\O_{\CPt}\oplus\O_{\CPt}(\Pt))$, we see immediately that $W^{\Pt}$ is naturally isomorphic 
to the surface $\Zt^{\Pt}$. Since $\Et$ is the direct image under $\pit:\Zt^{\Pt}\ra \CPt$ of $\Mt^{\Pt}(-\Pt)$ 
and $\tt$ is that of multiplication by $-\Xt^{\Pt}$, it follows also that
$\Xi^{\Pt}=(-1)_{rel}\widehat{\Sigma}^{\Pt}$. Finally, by the results of \cite{BNR} (Proposition 3.6 and Remark 3.7),
we obtain $Q^{\Pt}=(-1)_{rel}^{\ast}\Mt^{\Pt}(\Pt)$ as a sheaf. We know furthermore that the parabolic structure 
of $\Et$ induces a parabolic structure on $Q^{\Pt}$. On the other hand, because $\Mt^{\Pt}$ has a 
parabolic structure, the above isomorphism makes $Q^{\Pt}$ into a parabolic sheaf as well. These two parabolic
structures on $Q^{\Pt}$ agree: indeed, the parabolic structure of $\Et$ is the direct image of the one of 
$\Mt^{\Pt}$, so the filtration comes from the restriction of $\Mt^{\Pt}$ to some branches 
of the spectral curve over the parabolic points of $\Et$, hence the two filtrations of $Q^{\Pt}$ are the same; 
a similar argument works for the weights. 


The next ingredient in the construction is the analog of diagram (\ref{absdiag}). The surface $Z$ was obtained from 
$Z^\sP$ by blowing up the points $t^{+}$ of the $0$-section mapping to $\sP$ under $\pi$. Therefore, its analog $W$ for 
$W^{\Pt}$ is blow-up in the points $T^{+}$ of the $0$-section of $\pit$ over $\Pt$: clearly, this is the surface
$\Zt$. Now, $Z^{\abs}$ was obtained from $Z^\sP$ by blowing up the points $\Tt^{-}$ in the intersection of the
infinity-fiber of $\pi$ and the spectral curve. Because of $\Xi^{\Pt}=(-1)_{rel}\widehat{\Sigma}^{\Pt}$, the
intersection points of the infinity-fiber of $\pit$ and the spectral curve $\Xi^{\Pt}$ are the points
$(-1)^{\ast}T^{-}$ of $\Zt$. It follows that $(-1)$ induces a natural isomorphism between the absolute surface
$W^{\abs}$ of $(\Et,\tt)$ and $Z^{\abs}$; hence, we will simply write $W^{\abs}=(-1)_{rel}^{\ast}Z^{\abs}$. 
Notice that this surface has a projection to both projective lines $\CP$ and $\CPt$ (although 
these projections are only rational, and not every fiber is a single line), and the map $(-1)_{rel}$ above is 
induced by inversion on the fibers when it is considered as a fibration over $\CPt$. However, it is possible 
to interpret the same map as induced by inversion of the basis of the other fibration; we will simply write 
$(-1)$ for this map in the sequel, for any fibration over $\CP$. Therefore, $W^{\abs}$ can equally be written
as $(-1)^{\ast}Z^{\abs}$. We now come to an analog of $\Zt$: this surface was the blow-down in $Z^{\abs}$ of the 
proper transforms $\sE^{-}$ of the fibers of $\pi$ over the points $\sP$. Applying this to $W^{\abs}$ we obtain the result
$\Wt=(-1)^{\ast}Z$. Finally, arguments similar to the above yield that the analog of $\Zt^{\Pt}$ for $(\Et,\tt)$ is 
$\Wt^\sP=(-1)^{\ast}Z^\sP$, that is the surface ${\mathbf P}(\O_{\CP}\oplus\O_{\CP}((-1)^{\ast}\sP))$. We deduce that 
the diagram (\ref{absdiag}) corresponding to $(\Et,\tt)$ is 
\begin{align}\label{tabsdiag}
       \xymatrix{
       & & (-1)^{\ast}Z^{\abs} \ar[dl]_{(-1)^{\ast}\qmapPt} 
       \ar[dr]^{(-1)^{\ast}\qmapP} & & \\
       & \Zt \ar[dl]_{\pmaptP} \ar[dr]^{\pmaptz} & & 
       (-1)^{\ast}Z \ar[dl]_{(-1)^{\ast}\pmapz} \ar[dr]^{(-1)^{\ast}\pmapP} & \\
       \Zt^{\Pt} & & \CP \times \CPt & & (-1)^{\ast}Z^\sP
       }
\end{align}
The surface in the lower-left corner is the standard spectral surface of $(\Et,\tt)$, it has a projection to
$\CPt$, and the transformed bundle is obtained by taking the proper transform of $Q^{\Pt}$ in $(-1)^{\ast}Z^{\abs}$ with 
respect to $\pmaptP\circ (-1)^{\ast}\qmapPt$, 
deleting the exceptional divisor of $\pmaptP$ from the parabolic
divisor, adding the exceptional divisor of $(-1)^{\ast}\pmapP$ to the parabolic divisor, pushing down the
result to $(-1)^{\ast}Z^\sP$, then pushing down the result to $\CP$ by the projection map $\pi$ of $(-1)^{\ast}Z^\sP$, 
and finally tensoring by $(-1)^{\ast}\sP$. 

We have seen that the sheaf $Q^{\Pt}$ is isomorphic to 
$(-1)_{rel}^{\ast}\Mt^{\Pt}(\Pt)$. It follows from the
property 
$\wtd{(\pmaptP\circ \qmapPt)} \circ (\pmaptP\circ \qmapPt)_{\ast}= \mbox{Id}$ 
for pure sheaves of dimension $1$ on a smooth surface (see Lemma \ref{lem:Purity}), that 
$\wtd{(\pmaptP\circ \qmapPt)}Q^{\Pt}=(-1)^{\ast}\add{\widehat{\sE}^{+}}N^{\abs}$. To obtain the absolute spectral 
sheaf of $(\Et,\tt)$ we need to delete from the parabolic divisor of $(-1)^{\ast}\add{\widehat{\sE}^{+}}N^{\abs}$ 
the exceptional divisor $\widehat{\sE}^{+}$ of the blow-up map $\pmaptP$. 
We deduce from Proposition \ref{prop:adddelinv} that the absolute spectral 
sheaf of $(\Et,\tt)$ is $(-1)^{\ast}N^{\abs}$ on $(-1)^{\ast}Z^{\abs}$. The next step in the construction is 
to add the exceptional divisor $(-1)^{\ast}\sE^{+}$ of $(-1)^{\ast}\pmapP$ to the parabolic divisor of $(-1)^{\ast}N^{\abs}$. 
By Proposition \ref{prop:condts}, Condition \ref{assn:Main} for $(\E,\theta)$ implies Condition \ref{assn:second} for $N^{\abs}$ 
and $\sE^+$. Again by Proposition \ref{prop:adddelinv}, addition and deletion of a divisor are inverses to each other under 
Condition \ref{assn:second}. We obtain that 
$$
   \add{(-1)^{\ast}\sE^{+}}(-1)^{\ast}N^{\abs}=(-1)^{\ast}(\add{\sE^{+}}N^{\abs})=(-1)^{\ast}\wtd{(\pmaptP\circ \qmapPt)}M^P. 
$$ 
We then consider the direct image of this parabolic sheaf with respect to the blow-up map 
$(-1)^{\ast}(\pmaptP\circ \qmapPt)$: by Lemma \ref{lem:Purity}, the direct image is 
$(-1)^{\ast}M^\sP$. The push-down of this to $\CP$ by the projection $(-1)^{\ast}\pi$ of $(-1)^{\ast}Z^\sP$ 
is $(-1)^{\ast}(\E(\sP))=((-1)^{\ast}\E)((-1)^{\ast}\sP)$, see (\ref{lem:(-1)}). The final step is to tensor this sheaf by the 
inverse (in the divisor group) of the effective divisor corresponding to the parabolic set. Here this 
effective divisor is $(-1)^{\ast}\sP$, therefore tensoring $((-1)^{\ast}\E)((-1)^{\ast}\sP)$ by its inverse, we get 
precisely $(-1)^{\ast}\E$. This proves equality of the bundles $\Ett$ and $(-1)^{\ast}\E$. 
Clearly, the modifications of the sheaves involved so far transform the parabolic structure of $\Mt^{\Pt}$ 
into the parabolic structure of $(-1)^{\ast}M^\sP$ induced via pull-back by $(-1)$ of the original structure of 
$M^\sP$. Since the direct image by $\pi$ of the parabolic structure of $M^\sP(-\sP)$ is the parabolic structure 
of $\E$, we also see that the direct image by $\pi$ of the parabolic structure of $(-1)^{\ast}(M^\sP(-\sP))$ is the 
parabolic structure of $(-1)^{\ast}\E$ induced by pull-back under $(-1)$ from the parabolic structure of $\E$. 
Finally, the canonical section $\widehat{\widehat{x}}_{(-1)^{\ast}\sP}$ of $(-1)^{\ast}Z^\sP$ is $(-1)^{\ast}x_\sP$, 
 where $x_\sP$ is the canonical section of $Z^\sP$. By definition, the double transformed Higgs field $\ttt$ is the 
direct image with respect to $\pi$ of multiplication by $-\widehat{\widehat{x}}_{(-1)^{\ast}\sP}$. 
On the other hand,
the Higgs field $\theta$ is equal to the direct image of 
multiplication by $x_\sP$. It follows that $-\ttt=(-1)^{\ast}\theta$. 

\end{proof}

\section{The Map on Moduli Spaces}\label{sec:moduli}
In this section, we show that Nahm transform is a hyper-K\"ahler isometry of 
moduli spaces (Corollary \ref{isomcor}). 

As it is shown in Theorem 0.2 of \cite{BiqBoa}, the moduli space of stable Higgs bundles of parabolic degree 
$0$ with fixed simultaneously diagonalizable polar parts of arbitrary order and fixed parabolic structures is 
a hyper-K\"ahler manifold: the two anticommuting complex structures $J$ and $I$ are by definition given by the 
local holomorphic variations of Higgs bundles and integrable connections respectively, and the Euclidean
metric on the moduli space is defined as $L^{2}$ inner product of the harmonic representatives of tangent
vectors. Let now $\Mod$ denote the moduli space of Higgs bundles on $\CP$ with logarithmic singularities 
in the points of $\sD$ with fixed equivalence class of polar parts (\ref{thetaj})-(\ref{bj}) and fixed parabolic
filtration (\ref{parj}) and weights $\alpha^{j}_{k}$, and with an irregular singularity of rank one at infinity
with fixed equivalence class of polar parts (\ref{thetainf})-(\ref{binf}) and fixed parabolic structure with 
weights $\alpha^{\infty}_{k}$, up to complex gauge transformations preserving the parabolic structures. 
\begin{lem}\label{lem:moddim}
The complex dimension of the Zariski tangent space of $\Mod$ in any point is 
$$
2r\rt +2-r-\rt -\sum_{j=1}^{n}rk(res(\theta,p_{j}))^{2}-\sum_{l=1}^{\nt}rk(res(\tt,\xi_{l}))^{2},
$$ 
where the last two sums are taken for all logarithmic singularities of $\theta$ and $\tt$ respectively. 
\end{lem}
\begin{rk}
This formula is in fact invariant under exchanging $r$ with $\rt$ and $\theta$ with $\tt$, as it should be because 
of invertibility of the transform. 
\end{rk}
\begin{proof}
The computation is done in \cite{by} for the case of parabolic Higgs bundles of rank $r$ with only logarithmic 
singularities on a curve of arbitrary genus $g$. In fact, the authors there fix the residues of the Higgs field to be block 
nilpotent, but the same proof works for any other fixed block-diagonal parts as well.
The result obtained there is 
\begin{equation}\label{dimformula}
     2(g-1)r^{2}+2+\sum_{j=1}^{n}2f_{p_{j}}, 
\end{equation}
where $2f_{p_{j}}$ is the dimension of the adjoint orbit of $res(\theta,p_{j})$ in $\gl=\mathfrak{gl}(r,\kk)$ 
(see also the count of the dimension of the moduli space of parabolic vector bundles in \cite{ms}). 
We can understand this as coming from excision: the term $2(g-1)r^{2}+2$ is the degree of $\Om^{1}\otimes End(\E)$ 
plus the constant $2$ coming from global endomorphisms of $\E$, and in the last sum we add up terms arising 
in a neighborhood of each of the singular points. Explicitly, because we only consider deformations of the 
Higgs field whose residues in any singular point can be taken into the initial residue by a holomorphic change
of basis, the residue of an infinitesimal deformation corresponding to a one-parameter family of such 
deformations has to be in the adjoint orbit of the residue of $\theta$ in $\gl$; the dimension of such choices 
for the residue in $p_{j}$ is by definition $2f_{p_{j}}$. In the case where irregular singularities occur, 
by the same excision argument we need to define the quantity $2f_{p}$ in the last sum of (\ref{dimformula}) 
as the dimension of the adjoint orbit of the polar part of the Higgs field in 
$\gl\otimes_{\kk}\kk[\varepsilon]/(\varepsilon^{n_{p}+1})$, where $n_{p}$ is the Poincar\'e rank of the singularity in $p$. 

In our case, the only irregular singularity is infinity, of Poincar\'e rank $1$; let us compute the dimension 
of its orbit in $\gl\otimes_{\kk}\kk[\varepsilon]/(\varepsilon^{2})$. 
For the sake of simplicity, we only do this in the special case $\nt=2$; the generalization to higher $\nt$ is
immediate. So we suppose that at infinity the Higgs field can be written in the block form 
\begin{eqnarray}\label{block}
\frac{1}{2}\begin{pmatrix}
                \Xi_{1} &0\\
                0& \Xi_{2} 
           \end{pmatrix}+
           \varepsilon\begin{pmatrix}
                \Lambda_{1} & 0\\
                0 & \Lambda_{2}
            \end{pmatrix},    
\end{eqnarray}
where $\Xi_{1}=\xi_{1}\mbox{Id}_{a}$ for some $1\leq a < r$, $\Xi_{2}=\xi_{2}\mbox{Id}_{r-a}$,
$\Lambda_{1}=\mbox{diag}(\lambda^{\infty}_{1},\ldots,\lambda^{\infty}_{a})$, and $\Lambda_{2}=\mbox{diag}(\lambda^{\infty}_{1+a},\ldots,\lambda^{\infty}_{r})$ 
(see (\ref{a})-(\ref{binf})). We also assume $\xi_{1}\neq \xi_{2}$, and that all the $\lambda^{\infty}_{1},\ldots,\lambda^{\infty}_{a}$ 
are nonzero and pairwise distinct, and the same condition for $\lambda^{\infty}_{1+a},\ldots,\lambda^{\infty}_{r}$. 
Then the stabilizer of the adjoint action is by definition the block matrices 
\begin{eqnarray}
           \begin{pmatrix}
                A & B\\
                C& D 
           \end{pmatrix}+
           \varepsilon\begin{pmatrix}
                \alpha & \beta \\
                \gamma & \delta
            \end{pmatrix}
\end{eqnarray}
which commute with (\ref{block}) modulo $\varepsilon^{2}$. It is straightforward to check that this 
holds if and only if $B=0,C=0,\beta=0,\gamma=0$ and $A$ and $\sD$ are diagonal; under these assumptions, $\alpha$ and $\delta$ 
can be arbitrary. Therefore, the dimension of the stabilizer of the polar form (\ref{block}) is 
$a+(r-a)+a^{2}+(r-a)^{2}$, and so the dimension of its orbit is $2r^{2}-r-a^{2}-(r-a)^{2}$. For general 
$\nt$, the same argument gives for this dimension $2r^{2}-r-\sum_{l=1}^{\nt}(a_{l}-a_{l-1})^{2}$. Now, 
as the transformed Higgs field $\tt$ has residue of rank $(a_{l}-a_{l-1})$ in $\xi_{l}$, this can be 
rewritten as $2r^{2}-r-\sum_{l=1}^{\nt}rk(res(\tt,\xi_{l}))^{2}$. 

Similarly, it is easy to check that under the assumptions made in (\ref{bj}), for all logarithmic 
singularity $p_{j}$ the formula 
\begin{align*}
   2f_{p_{j}}&=r^{2}-r_{j}^{2}-(r-r_{j}) \\
             &=r^{2}-(r-rk(res(\theta,p_{j})))^{2}-rk(res(\theta,p_{j}))\\
             &=2r\cdot rk(res(\theta,p_{j})) - rk(res(\theta,p_{j}))^{2} - rk(res(\theta,p_{j})) \\
             &=(2r-1)rk(res(\theta,p_{j})) - rk(res(\theta,p_{j}))^{2}
\end{align*}
holds, where $r-r_{j}$ is the rank of $res(\theta,p_{j})$. Plugging these into (\ref{dimformula}) and using 
(\ref{trrank}) one obtains the dimension of the Zariski tangent 
as claimed. 
\end{proof}
Similarly to $\Mod$, let $\Modt$ denote the moduli space of stable Higgs bundles of parabolic degree $0$ on 
$\CPt$ with logarithmic singularities in the points of $\Dt$ with fixed equivalence class of polar parts 
and parabolic structures induced by the transform from the corresponding structures of $(\E,\theta)$ at infinity, 
and with an irregular singularity of rank one with fixed equivalence class of polar parts and fixed parabolic 
structure induced by the transform from the corresponding structures of $(\E,\theta)$ in the points of $\sD$ -- 
as explained in Section \ref{nahmsect} --, again up to complex gauge transformations preserving the parabolic 
structures. 

\begin{lem}\label{stablem}
Let $(\E_{\bullet},\theta)$ be a Higgs bundle on $\CP$ satisfying Condition
\ref{assn:Main}. 
Then the parabolic degrees of $\E_{\bullet}$ and of its Nahm transform $\Et_{\bullet}$ are the same. 
Furthermore, if $(\E_{\bullet},\theta)$ is stable of degree $0$, then the same is true for 
$(\Et_{\bullet},\tt)$. 
\end{lem}
\begin{proof}
The claim on parabolic degrees follows from Grothendieck-Riemann-Roch, as explained in Section 4.7 of
\cite{Sz}. It is also possible to deduce it using the fact 
that under the Condition \ref{assn:Main} all operations involved in passing from $M^\sP_{\bullet}$ to 
$\Mt^{\Pt}_{\bullet}$ preserve the parabolic 
Euler-characteristic. Indeed, by the parabolic Riemann-Roch theorem (Proposition 
\ref{prop:parRR}) applied to the curve $\CP$, 
one has $\para\chi(\E_{\bullet}(\sP+\infty))=\para\deg(\E_{\bullet})+r$, or equivalently 
$\para\chi(\E_{\bullet}(P))=\para\deg(\E_{\bullet})$. Of course, a similar relation holds for $\Et_{\bullet}$ as well. 
Finally, we obtain the result using that $\para\chi(\E_{\bullet}(\sP))=\para\chi(M^\sP_{\bullet})$ because 
$\pi_*M^\sP_{\bullet}=\E(\sP)_{\bullet}$, and the analogous statement for $\Et_{\bullet}$. 

Suppose now $(\E^{\prime}_{\bullet},\theta^{\prime})$ is a parabolic Higgs subbundle of $(\E_{\bullet},\theta)$. 
By Remark \ref{rk:indparstr}, we may suppose that the parabolic structure of $\E^{\prime}_{\bullet}$ is the 
structure induced by $\E_{\bullet}$. 
By Lemma \ref{lem:preserveSecAssn}, the standard spectral sheaf $(M^{\prime})^\sP$ of $\E^{\prime}_{\bullet}$ and the divisor 
$\sE^+$ also satisfy Condition \ref{assn:second}, because $(M^{\prime})^\sP_{\bullet}$ is a parabolic subsheaf of $M^\sP_{\bullet}$ 
with the induced parabolic structure. By Lemma \ref{lem:Injectivity}, proper transform 
preserves injective maps of sheaves. The same thing holds clearly for direct image by a blow-up map because it is the inverse 
of proper transform, and for addition and deletion, since on the level of sheaves these latter are simply tensoring operations. 
We conclude that $(\Mt^{\prime})^{\Pt}_{\bullet}$ is a parabolic subsheaf of $\Mt^{\Pt}_{\bullet}$, hence 
$(\Et^{\prime}_{\bullet},\tt^{\prime})$ is a parabolic Higgs subbundle of $(\Et_{\bullet},\tt)$. 
On the other hand, by the first part of the Lemma, the parabolic degree of $\Et^{\prime}_{\bullet}$ is equal to that of  
$\E^{\prime}_{\bullet}$. In particular, the parabolic degree of $\Et^{\prime}_{\bullet}$ is positive if and only if 
the parabolic degree of $\Et_{\bullet}$ is positive. 
In different terms, if $\para{\deg}(\E_{\bullet})=0$, then $(\E^{\prime}_{\bullet},\theta^{\prime})$ is 
destabilizing for $(\E_{\bullet},\theta)$ if and only if $(\Et^{\prime}_{\bullet},\tt^{\prime})$ is destabilizing for 
$(\Et_{\bullet},\tt)$. This proves preservation of stability. 
\end{proof}

The lemma allows us to introduce the map 
\begin{equation}\label{nahmmap}
    \Nahm : \Mod \lra \Modt
\end{equation}
defined by mapping the gauge equivalence class of the Higgs bundle $(\E,\theta)$ to the gauge equivalence class 
of the Higgs bundle $(\Et,\tt)$; by an easy adaptation of Lemma 1 of \cite{Jarsur} to the parabolic case over 
a curve, this map is well-defined. It is a bijective map between hyper-K\"ahler manifolds. We have the following
result. 
\begin{cor}\label{isomcor}
The map $\Nahm$ is a hyper-K\"ahler isometry. 
\end{cor}
\begin{proof}
By Theorem \ref{invthm}, the map $\Nahm$ is invertible. For the fact that $\Nahm$ preserves the
$L^{2}$-metric, we refer to \cite{BvB}: the computations there carry through to this case, 
because they only make use of the invertibility of the transform and general properties of the Green's 
operator that are satisfied in this case as well.

Therefore, all that remains is to check that it preserves the complex structures $I$ and $J$. 
Let us start with $J$: by Proposition 4.15 and equation (4.14) of \cite{Sz}, 
the restriction of the holomorphic bundle $\Et$ to the affine $\Ct$ is the first hypercohomology 
$\H^{1}(\E \xrightarrow{\theta_{\xi }} \F(\sP))$, together with the holomorphic structure induced by the trivial 
holomorphic structure of $\F(\sP)$ relative to $\CPt$. 
(Notice that what we called $\F$ in \cite{Sz} is called $\F(\sP)$ in the present paper.) 
Furthermore, by the extension \cite{Sz}, (4.35) of this holomorphic bundle to infinity, we have 
the following isomorphism of holomorphic bundles over $\CPt$:
\begin{equation}\label{hypcoh}
    \Et^{ind} = \H^{1}\left(\E \xrightarrow{y_{0}\theta -x_{0}} \F(\sP)\otimes\O_{\CPt}(1)\right), 
\end{equation}
where the right-hand side is endowed with the holomorphic structure induced from the holomorphic structure 
of the sheaf $\F(\sP)\otimes\O_{\CPt}(1)$ relative to $\CPt$. Let now $T$ be an open set in an affine complex line, 
and $(\E(t),\theta(t))$ with $t \in T$ be a local $1$-parameter holomorphic family of Higgs bundles with fixed 
singularity data. The transform maps each Higgs bundle in this family to a Higgs bundle $(\Et(t),\tt(t))$; 
we need to show that this is a holomorphic family of Higgs bundles on $\CPt$ over $T$. Because of 
equation (\ref{hypcoh}), the induced extensions of the $\Et$ vary holomorphically over $T$: indeed, 
the sheaves $\E(t)$ and $\F(t)$ depend holomorphically on $t$, as well as the map $y_{0}\theta(t) -x_{0}$, 
and the first hypercohomology spaces of a holomorphic family of sheaf complexes such that all the other 
hypercohomologies vanish, form again a holomorphic family. This implies 
that the bundles $\Et$ depend holomorphically on $t$ as well, because the parabolic structures are fixed, 
and in order to obtain the transformed extensions from the induced ones we only need to change the local 
holomorphic sections with non-zero weights, so this procedure does not depend on $t$ at all. 

On the other hand, by Theorem \ref{trthm} the map $\tt(t)$ is the direct image by $\pit(t)$ of multiplication by
$-\Xt^{\Pt}$; here $\pit(t)$ is the projection map from $\widehat{\Sigma}^{\Pt}(t)$ to $\CPt$. Now, since 
$\theta(t)$ varies holomorphically with $t$, it follows that so does the standard spectral curve $\Sigma^\sP(t)$ 
corresponding to $(\E(t),\theta(t))$. Since the eigenvalues of the residues of $\theta(t)$ are fixed, the elementary
transformations of the construction are performed at points that are independent of $t$. We deduce that 
the spectral curve $\widehat{\Sigma}^{\Pt}(t)$ and in particular the projection $\pit(t)$ depend holomorphically 
on $t$ as well. Since $\tt(t)$ is the direct image of multiplication by a coordinate that does not depend on 
$t$, with respect to a projection depending holomorphically on $t$, the resulting map depends also 
holomorphically on $t$. This proves that $\Nahm$ preserves the complex structure $J$. 

We now come to the case of the complex structure $I$: it can be treated in a very much similar way as $J$. 
Namely, in \cite{Sz2} the second author proves that the Nahm transform of a holomorphic bundle with integrable
connection $(E,\nabla)$ can also be given a holomorphic interpretation: by Proposition 5.1 of \emph{loc. cit.},
the holomorphic bundle $\widehat{E}$ on the open affine $\Ct$ is the first hypercohomology 
$\H^{1}(E\xrightarrow{\nabla_{\xi}}F)$, where $\nabla_{\xi}= \nabla -\xi \d z$ and $F$ is the sheaf generated by $E$ and 
the image of $\nabla$; and by Proposition 5.3 of \emph{loc. cit.}, the transformed integrable connection 
$\widehat{\nabla}$ is induced by the connection $\hat{\d}-z\d \xi$ relative to $\CPt$ on $F$. The rest of 
the argument follows the line of the case of the complex structure $J$. 
\end{proof}



\bibliographystyle{alpha}
\bibliography{spectral4}

\begin{thebibliography}{BNR89}

\bibitem[BB04]{BiqBoa}
Olivier Biquard and Philip Boalch.
\newblock Wild non-abelian {H}odge theory on curves.
\newblock {\em Compos. Math.}, 140(1):179--204, 2004.

\bibitem[BNR89]{BNR}
Arnaud Beauville, M.~S. Narasimhan, and S.~Ramanan.
\newblock Spectral curves and the generalised theta divisor.
\newblock {\em J. Reine Angew. Math.}, 398:169--179, 1989.

\bibitem[BvB89]{BvB}
Peter~J. Braam and Pierre van Baal.
\newblock Nahm's transformation for instantons.
\newblock {\em Comm. Math. Phys.}, 122(2):267--280, 1989.

\bibitem[BY96]{by}
Hans~U. Boden and K{\^o}ji Yokogawa.
\newblock Moduli spaces of parabolic {H}iggs bundles and parabolic {$K(D)$}
  pairs over smooth curves. {I}.
\newblock {\em Internat. J. Math.}, 7(5):573--598, 1996.

\bibitem[HL97]{HL}
Daniel Huybrechts and Manfred Lehn.
\newblock {\em The geometry of moduli spaces of sheaves}.
\newblock Aspects of Mathematics, E31. Friedr. Vieweg \& Sohn, Braunschweig,
  1997.

\bibitem[Jar04]{Jarsur}
Marcos Jardim.
\newblock A survey on {N}ahm transform.
\newblock {\em J. Geom. Phys.}, 52(3):313--327, 2004.

\bibitem[MS80]{ms}
V.~B. Mehta and C.~S. Seshadri.
\newblock Moduli of vector bundles on curves with parabolic structures.
\newblock {\em Math. Ann.}, 248(3):205--239, 1980.

\bibitem[Sza]{Sz2}
Szil\'{a}rd Szab\'{o}.
\newblock Transform\'ees de nahm et de laplace parabolique.
\newblock submitted.

\bibitem[Sza05]{Sz}
Szil\'{a}rd Szab\'{o}.
\newblock {\em Nahm transform for integrable connections on the Riemann
  sphere}.
\newblock PhD thesis, IRMA Strasbourg, 2005.

\bibitem[Yok93]{yo}
K{\^o}ji Yokogawa.
\newblock Compactification of moduli of parabolic sheaves and moduli of
  parabolic {H}iggs sheaves.
\newblock {\em J. Math. Kyoto Univ.}, 33(2):451--504, 1993.

\end{thebibliography}

\end{document}